\documentclass[journal]{IEEEtran}
\usepackage{amssymb}
\usepackage[pdftex]{graphicx}
\usepackage{tikz}

\usepackage{framed}
\usepackage{bm}
\usepackage{cite}
\usepackage{comment}
\usepackage{algorithmic}
\usepackage{algorithm}
\usepackage{enumerate}
\usepackage{enumitem}
\usepackage{amsmath}
\usepackage{amsfonts}
\usepackage{color}
\usepackage{cases}

\usetikzlibrary{intersections,calc,arrows.meta}

\newtheorem{theorem}{Theorem}
\newtheorem{assumption}{Assumption}

\newtheorem{remark}{Remark}
\newtheorem{lemma}{Lemma}
\newtheorem{corollary}{Corollary}

\newtheorem{proposition}{Proposition}

\newtheorem{problem}{Problem}

\newtheorem{example}{Example}
\newtheorem{proof}{Proof}

\usepackage{blkarray}

\usepackage[sort]{cleveref}
\usepackage{mathtools}
\usepackage{subcaption}
\usepackage{multirow}
\usepackage{threeparttable}

\usepackage{pgfplots}
\pgfplotsset{compat=1.18}
\usepgfplotslibrary{statistics}

\crefname{section}{Section}{Sections}
\crefname{appendix}{Appendix}{Appendices}
\crefname{theorem}{Theorem}{Theorems}
\crefname{proposition}{Proposition}{Propositions}
\crefname{lemma}{Lemma}{Lemmas}
\crefname{problem}{Problem}{Problems}
\crefname{assumption}{Assumption}{Assumptions}
\crefname{algorithm}{Algorithm}{Algorithms}
\crefname{remark}{Remark}{Remarks}
\crefname{equation}{}{}
\Crefname{equation}{Equation}{Equations}
\crefname{example}{Example}{Examples}
\crefname{figure}{Fig.}{Figs.}
\crefname{table}{Table}{Tables}

\crefname{enumi}{}{}

\newcommand{\RR}{\mathbb{R}}
\newcommand{\CC}{\mathbb{C}}
\newcommand{\NN}{\mathbb{N}}
\newcommand{\ZZ}{\mathbb{Z}}
\newcommand{\norm}[1]{\left\Vert{#1}\right\Vert}
\newcommand{\rd}{\mathrm{d}}
\newcommand{\tblue}[1]{\textcolor{blue}{#1}}

\newcommand{\e}{\mathrm{e}}

\DeclareMathOperator*{\minimize}{minimize}
\DeclareMathOperator*{\argmin}{arg\,min}
\DeclareMathOperator*{\dom}{dom}
\DeclareMathOperator*{\interior}{int}
\DeclareMathOperator{\diag}{diag}
\DeclareMathOperator{\rank}{rank}
\DeclareMathOperator{\vect}{vec}
\DeclareMathOperator{\tr}{tr}
\DeclareMathOperator{\re}{Re}
\DeclareMathOperator{\im}{Im}

\crefname{ALC@unique}{Step}{Steps}

\usepackage{latexsym}
\def\qed{\hfill $\Box$} 

\usepackage{cite}
\usepackage{amsmath}
\usepackage{url}

\begin{document}
%
\title{Infinite-horizon controllability scores for\\ linear time-invariant systems}
\author{Kota Umezu and Kazuhiro Sato\thanks{K. Umezu and K. Sato are with the Department of Mathematical Informatics, Graduate School of Information Science and Technology, The University of Tokyo, Tokyo 113-8656, Japan, email: krr0814bz@g.ecc.u-tokyo.ac.jp (K. Umezu), kazuhiro@mist.i.u-tokyo.ac.jp (K. Sato) }}
\maketitle
\thispagestyle{empty}
\pagestyle{empty}

\begin{abstract}
We formulate infinite-horizon controllability scores for linear
time-invariant systems whose conventional controllability Gramians may
diverge.
A mode-dependent scaling yields finite-horizon volumetric
controllability score (VCS) and average energy controllability score
(AECS) problems that are exactly equivalent to the original ones,
while the scaled Gramians converge as the horizon tends to infinity.
Using epi-convergence, we establish convergence of the finite-horizon
optimal allocations to the limiting solution set and obtain convergence
to a unique limiting score under explicit rank conditions.
The limiting VCS retains the stable, center, and unstable
components and preserves full-system controllability, whereas the
limiting AECS depends only on the stable component and may neglect an
input required for a center or unstable mode.
The VCS theory allows arbitrary spectral structures, whereas the AECS limit requires a nonempty stable component. The additional restriction on the center spectrum is imposed only for the real-Schur implementation.
A directed-Laplacian example verifies the convergence and uniqueness
conditions and illustrates the controllability distinction between the
two limiting scores.
\end{abstract}

\begin{IEEEkeywords}
controllability score, infinite horizon, linear time-invariant system, network centrality
\end{IEEEkeywords}

\IEEEpeerreviewmaketitle

\section{Introduction}
\label{sec:introduction}

Large-scale networks play an essential role in understanding complex phenomena and in intervening in complex dynamical systems across diverse domains~\cite{Newman2003}.
Many systems of contemporary interest---such as brain networks~\cite{BresslerMenon2010,WigSchlaggarPetersen2011}, biological interaction networks~\cite{BarabasiGulbahceLoscalzo2011}, social networks~\cite{Backstrom2012}, and power grids~\cite{Buldyrev2010,PaganiAiello2013}---are represented as networks consisting of a large number of interacting components.
In such systems, macroscopic behavior frequently emerges from collective interactions rather than from individual units alone.
While structural descriptions are informative, they are frequently insufficient to predict or steer real-world network dynamics, motivating frameworks that explicitly incorporate dynamics in addition to topology~\cite{Liu2011,Liu2012,ElterenQuaxSloot2022}.

Control theory, and in particular controllability, provides a dynamics-aware framework for analyzing and intervening in complex networked systems~\cite{Liu2011,Liu2012}. Controllability characterizes whether a system can be steered to desired states~\cite{Kalman1960} and leads to quantitative descriptors that reflect both dynamics and interconnection structure. Controllability-based analyses have been shown to be useful in
neuroscience \cite{Gu2015} and biology \cite{Vinayagam2016}, among
others. These developments motivate node-level measures that quantify where to intervene and how strongly, while accounting for the control effort required to influence the network.

In this context, controllability scores~\cite{SatoTerasaki2024} provide controllability-based centrality measures for continuous-time linear time-invariant (LTI) dynamical networks of the form
\begin{equation}\label{eq:lti_base}
    \dot{x}(t)=Ax(t),
\end{equation}
where $x(t)\in\RR^n$ is the state vector collecting the states of $n$ nodes, and the system matrix $A=(a_{ij})\in\RR^{n\times n}$ captures weighted interactions:
$a_{ij}$ represents the influence of node $j$ on node $i$.
We use ``node $i$" and ``state $x_i$" interchangeably.
Under mild conditions, controllability scores are defined as the unique optimal solutions of convex optimization problems, thereby assigning each node a quantitative measure of its importance for controlling the network dynamics.

The controllability scores address limitations of existing controllability-based frameworks, which can be broadly categorized into structural controllability~\cite{Lin1974} and quantitative controllability measures~\cite{MullerWeber1972}. Structural controllability depends only on network topology and enables graph-theoretic analysis without identifying edge weights~\cite{Liu2011,Liu2012}, but it does not preclude prohibitively large control-energy requirements~\cite{Yan2015}. Quantitative measures incorporate edge weights and time horizons to capture energy-related aspects of controllability; in this spirit, controllability scores quantify node importance while accounting for control effort. Existing Gramian-based centralities such as VCE and ACE introduced in \cite{Summers2016} may still be misaligned with quantitative controllability because they rely on the rank or pseudoinverse of a single-input controllability Gramian~\cite{SatoTerasaki2024}. By contrast, controllability scores are defined by introducing virtual actuation at all state nodes and optimizing the input allocation.
Subsequent studies established uniqueness properties~\cite{SatoTerasaki2024,SatoKawamura2025}, demonstrated applicability to human brain networks~\cite{SatoKawamura2025}, and incorporated constraints on actuated nodes~\cite{Sato2025}.

In \cref{subsec:controllability_score}, we define two controllability scores, VCS and AECS, as optimal solutions to distinct controllability-Gramian--based optimization problems
over a finite horizon $T>0$.
When a control deadline is prescribed, the corresponding finite-horizon score is appropriate. In the absence of a natural deadline, however, a horizon-independent reference is useful for assessing long-term actuator importance.
Previous work has shown that short horizons emphasize local,
short-term effects, whereas the allocations approach
network-dependent limiting patterns as $T$ increases
\cite{SatoKawamura2025}.
For a non-Hurwitz system, however, the passage to $T=\infty$ is
singular.
Stable, center (imaginary-axis), and unstable spectral components have
different asymptotic growth rates, and the conventional controllability
Gramian may grow polynomially or exponentially.
Consequently, its log determinant, inverse, and associated gradients
can become numerically unreliable even before overflow occurs.
The infinite-horizon problems therefore cannot be obtained by simply
substituting $T=\infty$ into the finite-horizon formulas.
This paper addresses this issue.

The contributions of this paper are summarized as follows:
\begin{enumerate}
    \item
    We introduce a mode-dependent scaling that accommodates stable,
    imaginary-axis, and unstable Jordan blocks, and prove convergence
    of the resulting scaled controllability Gramian.
    For every finite $T>0$, the scaled VCS and AECS problems are
    exactly equivalent to the original problems.
    Thus, the scaling is not an approximation and does not alter their
    optimal allocations.

    \item
    We formulate convex optimization problems whose optimal solutions define the limiting VCS and AECS, and establish epi-convergence of the corresponding extended-real-valued finite-horizon objectives.
    Every cluster point of the finite-horizon scores is therefore an optimal solution of the limiting problem, and the distance to the limiting solution set converges to zero. When the limiting score is unique, the finite-horizon score converges to it as $T\to\infty$.
   This variational argument accommodates the strict positive-definiteness constraints and possible changes in the effective domain in the infinite-horizon limit.

    \item
    We derive full-column-rank conditions for
    uniqueness of the limiting scores and identify structural
    mechanisms by which these conditions can fail.
    Uniqueness is not inherited automatically from the finite-horizon
    problems: we provide an example in which every finite-horizon
    optimizer is unique, whereas the limiting problem has the entire
    simplex as its solution set.

    \item
    We show that the limiting VCS retains the stable, center,
    and unstable components and thereby preserves full-system
    controllability through its feasibility condition.
    In contrast, the limiting AECS depends only on the stable
    component.
    It may therefore assign zero weight to an input that is unnecessary
    for controlling the stable modes but indispensable for a center or
    unstable direction.

    \item
    For systems with arbitrary stable and unstable eigenvalues and no imaginary-axis eigenvalues except possibly a semisimple zero, we replace the analytical Jordan representation by a computational procedure based on real Schur decompositions and Sylvester and Lyapunov equations.
    This avoids the numerically ill-conditioned computation of the Jordan canonical form and enables the limiting problems to be implemented using standard numerical linear-algebra routines.
\end{enumerate}

The remainder of this paper is organized as follows:
\Cref{sec:preliminaries} reviews finite-horizon controllability
scoring and identifies the challenges associated with the
infinite-horizon limit.
\Cref{sec:reformulation,sec:extension} develop the scaled formulation
and the limiting problems, together with their convergence and
uniqueness properties.
\Cref{sec:algorithms} presents the real-Schur implementation and its
computational scope, and \cref{sec:experiments} provides numerical
validation.
Finally, \cref{sec:conclusion} concludes the paper.

\subsubsection*{Notation}
Throughout this paper, $n$ denotes the number of nodes in the network, and $A\in\RR^{n\times n}$ denotes the system matrix in \cref{eq:lti_base}.
The sets of all natural numbers, all real numbers, and all complex numbers are denoted by $\NN$, $\RR$, and $\CC$, respectively.
The symbol $\sqrt{-1}$ denotes the imaginary unit.
Let $e_i\coloneq(0,\ldots,0,1,0,\ldots,0)^\top\in\RR^n$ be a standard vector whose $i$-th entry is $1$ and other entries are $0$.
Let $I_k$ and $N_k$ denote the identity matrix of size $k$ and the nilpotent Jordan block of size $k$ with ones on the superdiagonal, respectively.
Let $O$ denote a zero matrix of a suitable size.
For $v_1,\ldots,v_\ell\in\CC$,
$\mathrm{diag}(v_1,\ldots,v_\ell)$ denotes the diagonal matrix with
diagonal entries $v_1,\ldots,v_\ell$.
More generally, for square matrices
$X_i\in\CC^{k_i\times k_i}$ $(i=1,\ldots,m)$,
$\mathrm{diag}(X_1,\ldots,X_m)$ denotes the block-diagonal matrix
whose diagonal blocks are $X_1,\ldots,X_m$.
Let $\Delta$ be the standard simplex in $\RR^n$, i.e.,
\begin{equation*}
    \Delta\coloneq
    \left\{
        (p_i)\in\RR^n
    \,\middle|\,
        \begin{gathered}
            p_i\geq0 \ (i=1,\ldots,n), \\
            \textstyle \sum_{i=1}^n p_i=1
        \end{gathered}
    \right\}.
\end{equation*}

For $\lambda\in\CC$, $\re(\lambda)$, $\im(\lambda)$, and $\overline{\lambda}$ denote the real part, the imaginary part, and the complex conjugate of $\lambda$, respectively.
For a real matrix $X\in\RR^{k\times\ell}$, $X^\top$ denotes the transpose of $X$, and for a complex matrix $Y\in\CC^{k\times\ell}$, $Y^*$ denotes the conjugate transpose of $Y$.
For a real invertible matrix $X\in\RR^{\ell\times\ell}$, we write $X^{-\top}\coloneq(X^{-1})^\top$, and for a complex invertible matrix $Y\in\CC^{\ell\times\ell}$, we write $Y^{-*}\coloneq(Y^{-1})^*$.
For $Z\in\CC^{\ell\times\ell}$, $\e^Z, \det Z$, and $\tr Z$ denote the exponential, the determinant, and the trace of $Z$, respectively.
For $C=(c_{ij})\in\CC^{k\times\ell}$ and $C'\in\CC^{k'\times \ell'}$, $C\otimes C'\coloneq\begin{pmatrix} c_{11} C' & \cdots & c_{1\ell} C' \\ \vdots & \ddots & \vdots \\ c_{k1} C' & \cdots & c_{k\ell} C' \end{pmatrix} \in\CC^{kk'\times\ell\ell'}$ denotes their Kronecker product.
For Hermitian matrices $X,Y\in\CC^{\ell\times\ell}$, we write $X\succ Y$ if $X-Y$ is positive definite.
For $v=(v_1,\ldots,v_\ell)^\top\in\CC^\ell$, $\norm{v}\coloneq\sqrt{\sum_{i=1}^\ell |v_i|^2}$ denotes the standard Euclidean norm, and for $X\in\CC^{k\times\ell}$, $\norm{X}\coloneq\sqrt{\tr(X^* X)}$ denotes the Frobenius norm.

\section{Preliminaries}
\label{sec:preliminaries}

\subsection{Controllability score}
\label{subsec:controllability_score}

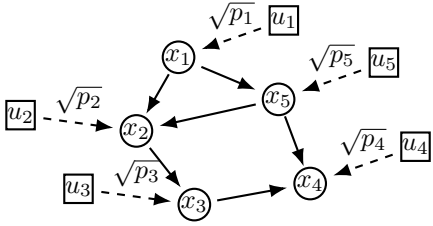
\begin{figure}[htbp]
    \centering
    \begin{tikzpicture}[scale=0.7]
        \node (P1) at (-0.3,2.8) {$x_1$};
        \node (P2) at (-1.1,1.4) {$x_2$};
        \node (P3) at (0,0) {$x_3$};
        \node (P4) at (2.2,0.4) {$x_4$};
        \node (P5) at (1.6,2) {$x_5$};

        \coordinate (UR) at (2,0.7);
        \coordinate (UL) at (-2.2,0.3);

        \node (U1) at ($(P1)+(UR)$) {$u_1$};
        \node (U2) at ($(P2)+(UL)$) {$u_2$};
        \node (U3) at ($(P3)+(UL)$) {$u_3$};
        \node (U4) at ($(P4)+(UR)$) {$u_4$};
        \node (U5) at ($(P5)+(UR)$) {$u_5$};

        \draw[thick] (P1) circle[radius=0.30];
        \draw[thick] (P2) circle[radius=0.30];
        \draw[thick] (P3) circle[radius=0.30];
        \draw[thick] (P4) circle[radius=0.30];
        \draw[thick] (P5) circle[radius=0.30];

        \draw[thick] ($(U1) + (-0.27,-0.27)$) rectangle +(0.54,0.54);
        \draw[thick] ($(U2) + (-0.27,-0.27)$) rectangle +(0.54,0.54);
        \draw[thick] ($(U3) + (-0.27,-0.27)$) rectangle +(0.54,0.54);
        \draw[thick] ($(U4) + (-0.27,-0.27)$) rectangle +(0.54,0.54);
        \draw[thick] ($(U5) + (-0.27,-0.27)$) rectangle +(0.54,0.54);

        \draw[-Latex,thick] (P1) -- (P2);
        \draw[-Latex,thick] (P1) -- (P5);
        \draw[-Latex,thick] (P2) -- (P3);
        \draw[-Latex,thick] (P3) -- (P4);
        \draw[-Latex,thick] (P5) -- (P2);
        \draw[-Latex,thick] (P5) -- (P4);

        \draw[-Latex,thick,dashed] (U1) -- (P1) node[midway, above] {$\sqrt{p_1}$};
        \draw[-Latex,thick,dashed] (U2) -- (P2) node[midway, above] {$\sqrt{p_2}$};
        \draw[-Latex,thick,dashed] (U3) -- (P3) node[midway, above] {$\sqrt{p_3}$};
        \draw[-Latex,thick,dashed] (U4) -- (P4) node[midway, above] {$\sqrt{p_4}$};
        \draw[-Latex,thick,dashed] (U5) -- (P5) node[midway, above] {$\sqrt{p_5}$};
    \end{tikzpicture}
    \caption{Virtual system for defining controllability scores.
    Each node $x_i$ is paired one-to-one with a virtual input node $u_i$; the allocation $p$ distributes a fixed control budget, and $u_i$ acts on $x_i$ with gain $\sqrt{p_i}$.}
    \label{fig:virtual_system}
\end{figure}

We begin with an intuitive description of controllability scores~\cite{SatoTerasaki2024} to motivate the formal definitions introduced later in this subsection.
To this end, we consider
the autonomous system \eqref{eq:lti_base} and
introduce a virtual actuation mechanism (see \cref{fig:virtual_system}).
Specifically, each state node $x_i$ is paired one-to-one with a virtual input $u_i$, and the input strength is determined by an allocation vector $p\in\Delta$:
\begin{equation}\label{eq:lti_diagonal}
    \dot{x}(t)=Ax(t)+\diag(\sqrt{p_1},\ldots,\sqrt{p_n})u(t),
\end{equation}
where $u(t)\coloneq(u_1(t),\ldots,u_n(t))^\top\in\RR^n$.
The value $p_i$ represents the amount of control resources allocated to node $x_i$; under a constraint on the input energy, a larger value of $p_i$ allows a stronger control effect on node $x_i$, indicating that node $x_i$ is a central node in the system.
The simplex constraint has a direct control-resource interpretation.
The value $p_i$ represents the fraction of a fixed total squared
input-gain budget assigned to node $i$.
The conditions $p_i\geq0$ ensure real nonnegative actuator gains,
whereas the equality constraint removes arbitrary common scaling.
Accordingly, projection onto $\Delta$ restores a trial allocation to
the admissible nonnegative fixed-budget set; it does not project the
physical state or input.
Given a controllability measure for system \cref{eq:lti_diagonal}, we optimize it over $p$ and denote an optimal allocation by $p^*$.
This optimal allocation $p^*$ maximizes controllability, and its entries $p_i^*$ quantify node importance.
We therefore interpret $p^*$ as a dynamics-aware centrality measure, called a controllability score.

To turn this intuition into a formal definition, we introduce controllability-Gramian--based scalar measures of controllability as functions of $p$.
We first define, for $B\in\CC^{n\times n}$ and $T>0$,
\begin{equation}
    \widetilde{W}(B;T)\coloneq\int_{0}^{T} \e^{At}B\e^{A^\top t}\rd t.
\end{equation}
The controllability Gramian for system \cref{eq:lti_diagonal} with terminal time $T>0$ is given by $\widetilde{W}(\diag(p);T)$, which we abbreviate as $\widetilde{W}(p;T)$.
For later algebraic arguments, the notation $\widetilde{W}(p;T)$ is
extended linearly to $p\in\RR^n$, without associating a physical input
matrix with vectors having negative components.

The quantities $-\log\det \widetilde{W}(p;T)$ and $\tr\left(\widetilde{W}(p;T)^{-1}\right)$ can serve as scalar controllability measures.
A smaller value of $-\log\det \widetilde{W}(p;T)$ corresponds to a larger volume of the reachable set $\mathcal{E}(p;T)$, defined as
\begin{equation*}
    \mathcal{E}(p;T)\coloneq
    \left\{
        x(T)\in\RR^n
    \,\middle|\,
        \begin{gathered}
            \textstyle\int_0^T\norm{u(t)}^2\rd t\leq1,\ x(0)=0 \\
            \text{under \cref{eq:lti_diagonal}}
        \end{gathered}
    \right\}.
\end{equation*}
A smaller value of $\tr\left(\widetilde{W}(p;T)^{-1}\right)$ implies a smaller average minimum energy to steer the state from the origin to points on the unit sphere $\{y\in\RR^n\mid\norm{y}= 1\}$, where the average is taken with respect to the uniform distribution on the sphere.

We thus consider the following minimization problems, which we refer to as controllability scoring problems.
\begin{problem}\label{prob:vcs}
    \begin{align*}
        \minimize_{p} \quad & {-}\log\det \widetilde{W}(p;T) \\
        \mathrm{subject\ to} \quad & p\in\Delta,\ \widetilde{W}(p;T)\succ O.
    \end{align*}
\end{problem}

\begin{problem}\label{prob:aecs}
    \begin{align*}
        \minimize_{p} \quad & \tr\left(\widetilde{W}(p;T)^{-1}\right) \\
        \mathrm{subject\ to} \quad & p\in\Delta,\ \widetilde{W}(p;T)\succ O.
    \end{align*}
\end{problem}
The optimal solution to \cref{prob:vcs} is referred to as the volumetric controllability score (VCS), whereas that of \cref{prob:aecs} is referred to as the average energy controllability score (AECS).
As explained above, these quantities can be interpreted as measures of the importance of individual nodes from the perspective of controllability.

\subsection{Remaining challenges and our approach}
\label{subsec:challenges}

The finite- and infinite-horizon scores serve different purposes.
When a control deadline is prescribed, that value should be used in
\cref{prob:vcs,prob:aecs}.
When no natural deadline is available, however, the choice of $T$
is arbitrary and can impair the reproducibility of long-term actuator
evaluation.
The infinite-horizon allocation provides a horizon-independent
reference, but raises the following two challenges:

\begin{enumerate}
  \item\label{challenge:1}
    For a general non-Hurwitz $A$, the unscaled Gramian does not
    converge to a finite matrix as $T\to\infty$.
    Hence, simply setting $T=\infty$ in the finite-horizon formulations
    does not yield well-defined infinite-horizon optimization problems.
    Conventional nodewise Gramian centralities face the same difficulty,
    as detailed in \cref{rem:centrality_divergence}.

  \item\label{challenge:2}
    Even at a large but finite $T$, the unscaled Gramian may contain
    entries and eigenvalues on widely separated scales, causing overflow,
    loss of numerical positive definiteness, and inaccurate evaluations
    of determinants, inverses, and gradients.
\end{enumerate}

We address both issues by scaling the Gramian according to the
asymptotic behavior of each spectral block: imaginary-axis Jordan
blocks produce polynomial growth, whereas unstable blocks produce
exponential growth.
The resulting matrix $W(p;T)$ converges, while the associated
finite-horizon programs remain exactly equivalent to
\cref{prob:vcs,prob:aecs}.
We then formulate the limiting objectives as extended-real-valued
convex functions and establish epi-convergence.
This justifies the convergence properties of the finite-horizon VCS
and AECS even when the strict positive-definiteness constraints lead
to changes in the effective domains.

The Jordan form is used only to characterize the general
infinite-horizon limit.
The limiting VCS formulation imposes no spectral restriction, whereas
the limiting AECS requires at least one eigenvalue with negative real
part.
The additional assumption in \cref{sec:algorithms} is needed only for
the real-Schur implementation, which avoids computing the Jordan form.
Because the scaled and unscaled programs have the same set of optimal
solutions for every finite $T$, the scaled formulation also provides a
numerically safer representation at large horizons.
It is not claimed to be uniformly cheaper than direct computation based
on the unscaled Gramian; the computational trade-off is discussed in
\cref{sec:algorithms,sec:experiments}.

\begin{remark}\label{rem:centrality_divergence}
    The AC, VCE, and ACE proposed in \cite{Summers2016} are nodewise
    functions of the finite-horizon Gramians
    $\widetilde{W}(e_ie_i^\top;T)$.
    Their conventional infinite-horizon definitions require
        $\widetilde{W}_{i,\infty}
        \coloneq
        \int_0^\infty
        \e^{At}e_ie_i^\top\e^{A^\top t}\rd t$
    to be finite for every $i$.
    If $A$ is non-Hurwitz, there exists an index $i$ such that $e_i$
    has a nonzero component in the spectral subspace associated with
    eigenvalues satisfying $\re(\lambda)\geq0$.
    Thus, $\widetilde{W}_{i,\infty}$ is not finite, and these measures do not define an infinite-horizon centrality value for every node in the network.
    For this $i$, the finite-horizon AC and VCE diverge as
    $T\to\infty$.
    Although ACE may have a finite scalar limit, that limit is not
    obtained from a finite infinite-horizon Gramian and suppresses the
    contribution of the diverging directions.

    VCE and ACE also omit zero Gramian eigenvalues while penalizing
    arbitrarily small positive ones \cite{SatoTerasaki2024}.
    By contrast, VCS and AECS are optimal allocations over a bounded
    feasible set and can converge even when the unscaled Gramian
    diverges.
    Hence, \cref{sec:experiments} focuses on their finite-to-infinite
    convergence rather than assigning $T=\infty$ values to AC, VCE,
    or ACE.
\end{remark}

\section{Scaling method of the controllability Gramian and reformulation of \cref{prob:vcs,prob:aecs}}
\label{sec:reformulation}

To define the infinite-horizon controllability scores in \cref{sec:extension}, we first introduce a scaling of the controllability Gramian and derive numerically stable reformulations of the finite-horizon \cref{prob:vcs,prob:aecs}.
Specifically, in this section we reformulate \cref{prob:vcs} as
\begin{problem}\label{prob:scaled_vcs}
    \begin{align*}
        \minimize_{p} \quad & {-}\log\det W(p;T) \\
        \mathrm{subject\ to} \quad & p\in\Delta,\ W(p;T)\succ O,
    \end{align*}
\end{problem}
and \cref{prob:aecs} as
\begin{problem}\label{prob:scaled_aecs}
    \begin{align*}
        \minimize_{p} \quad & \tr\left(Q^{-*}D(T)^{-*}W(p;T)^{-1}D(T)^{-1}Q^{-1}\right) \\
        \mathrm{subject\ to} \quad & p\in\Delta,\ W(p;T)\succ O.
    \end{align*}
\end{problem}
Here, $W(p;T), Q$, and $D(T)$ are defined in \cref{subsec:settings}.
We also prove that $W(p;T)$ converges as $T\to\infty$, and discuss the advantages of \cref{prob:scaled_vcs,prob:scaled_aecs}.

\subsection{Problem settings}
\label{subsec:settings}

Throughout this section, the matrix $A\in\RR^{n\times n}$ is arbitrary, and we assume that it admits a Jordan decomposition
    $Q^{-1}AQ=J\coloneq
    \mathrm{diag}(J_1,\ldots,J_m)$
for some invertible matrix $Q\in\CC^{n\times n}$.
Here, $m$ is the number of Jordan blocks, and $J_i=\lambda_i I_{n_i}+N_{n_i}\in\CC^{n_i\times n_i}$ denotes a Jordan block of size $n_i$ associated with the eigenvalue $\lambda_i$ for $i=1,\ldots,m$.
Let $m_-, m_0$, and $m_+$ denote the number of Jordan blocks satisfying $\re(\lambda_i)<0, \re(\lambda_i)=0$, and $\re(\lambda_i)>0$, respectively.
Moreover, let the eigenvalues be ordered such that
\begin{equation*}
    \re(\lambda_i)
    \begin{cases*}
        <0 & if $i=1,\ldots,m_-$, \\
        =0 & if $i=m_-+1,\ldots,m_-+m_0$, \\
        >0 & if $i=m_-+m_0+1,\ldots,m$,
    \end{cases*}
\end{equation*}
and
\begin{equation*}
    \im(\lambda_{m_-+1})\leq\im(\lambda_{m_-+2})\leq\cdots\leq\im(\lambda_{m_-+m_0}).
\end{equation*}

Let $D_i(T)\in\CC^{n_i\times n_i}$ be defined by
\begin{equation*}
    D_i(T)\coloneq\left\{
    \begin{aligned}
        & I_{n_i} \qquad \text{if $i=1,\ldots,m_-$}, \\
        & \diag(T^{n_i-\frac{1}{2}},T^{n_i-\frac{3}{2}},\ldots,T^{\frac{1}{2}}) \\
        & \qquad \qquad \text{if $i=m_-+1,\ldots,m_-+m_0$}, \\
        & \e^{J_i T} \qquad \text{if $i=m_-+m_0+1,\ldots,m$},
    \end{aligned}
    \right.
\end{equation*}
and define a scaling matrix $D(T)\in\CC^{n\times n}$ by
\begin{equation}\label{eq:D_definition}
    D(T)\coloneq \mathrm{diag}(D_1(T),\ldots,D_m(T)).
\end{equation}
We then define the scaled matrix
\begin{equation}\label{eq:scaledW_definition}
    W(B;T)\coloneq D(T)^{-1}Q^{-1}\widetilde{W}(B;T)Q^{-*}D(T)^{-*}.
\end{equation}
In what follows, we abbreviate $W(\diag(p);T)$ as $W(p;T)$ for $p\in\RR^n$, and refer to it as the scaled controllability Gramian.

We further introduce additional notation for later use.
Let $\displaystyle n_-\coloneq\sum_{i=1}^{m_-} n_i,\ n_0\coloneq\sum_{i=m_-+1}^{m_-+m_0} n_i$, and $\displaystyle n_+\coloneq\sum_{i=m_-+m_0+1}^{m} n_i$ denote the numbers of eigenvalues with negative, zero, and positive real parts, respectively, counted with algebraic multiplicity.
Let
    $J_-\coloneq
    \mathrm{diag}(J_1,\ldots,J_{m_-})
    \in\CC^{n_-\times n_-}$
be the Jordan matrix corresponding to the eigenvalues with negative real parts.
Similarly, we define $J_0\in\CC^{n_0\times n_0}$ and $J_+\in\CC^{n_+\times n_+}$.

Furthermore, let $m'_0$ denote the number of distinct values among $\lambda_{m_-+1},\ldots,\lambda_{m_-+m_0}$, and let $\lambda'_1,\cdots,\lambda'_{m'_0}$ be these values, ordered so that $\im(\lambda'_1)<\cdots<\im(\lambda'_{m'_0})$.
We define the geometric multiplicity of $\lambda'_i$ as $\displaystyle m_0^{(i)}\coloneq\sum_{\lambda_j=\lambda'_i} 1$, that is, the number of Jordan blocks associated with $\lambda_i'$, and the algebraic multiplicity of $\lambda'_i$ as $\displaystyle n_0^{(i)}\coloneq\sum_{\lambda_j=\lambda'_i} n_j$.

\subsection{Infinite-horizon scaled controllability Gramian}
\label{subsec:scaled_gramian}

We discuss the limit of the scaled controllability Gramian $W(B;T)$ as $T\to\infty$.
We repeatedly use the following lemma.

\begin{lemma}[{\cite[Theorem 2.4.4.1]{HornJohnson2012}}]\label{prop:sylvester_uniqueness}
    Let $C_1\in\CC^{\ell_1\times\ell_1}, C_2\in\CC^{\ell_2\times\ell_2}$, and $P\in\CC^{\ell_1\times\ell_2}$.
    If $\mu_1+\mu_2\neq0$ holds for any eigenvalues $\mu_1$ of $C_1$ and $\mu_2$ of $C_2$, then
        $C_1 Y+YC_2+P=O$
    has a unique solution $Y\in\CC^{\ell_1\times\ell_2}$.
\end{lemma}

Let us define $\widetilde{B}\coloneq Q^{-1}BQ^{-*}$ and partition $\widetilde{B}$ into blocks in the following two ways:
\begin{equation}\label{eq:Btilde_partition}
    \widetilde{B}=
    \begin{pmatrix}
        \widetilde{B}_{1,1} & \cdots & \widetilde{B}_{1,m} \\
        \vdots & \ddots & \vdots \\
        \widetilde{B}_{m,1} & \cdots & \widetilde{B}_{m,m}
    \end{pmatrix} \\
    =
    \begin{pmatrix}
        \widetilde{B}_{-,-} & \widetilde{B}_{-,0} & \widetilde{B}_{-,+} \\
        \widetilde{B}_{0,-} & \widetilde{B}_{0,0} & \widetilde{B}_{0,+} \\
        \widetilde{B}_{+,-} & \widetilde{B}_{+,0} & \widetilde{B}_{+,+}
    \end{pmatrix},
\end{equation}
where $\widetilde{B}_{i,j}\in\CC^{n_i\times n_j} \ (i,j=1,\ldots,m)$, and $\widetilde{B}_{\alpha,\beta}\in \CC^{n_\alpha\times n_\beta} \ (\alpha,\beta\in\{-,0,+\})$.
Then, let $W_-(B)$ and $W_+(B)$ be the unique solutions to
\begin{subequations}\label{eq:W_jordan_sylvester}
\begin{align}
    J_-W_-(B)+W_-(B)J_-^*+\widetilde{B}_{-,-}&=O, \label{eq:W_stable_jordan_sylvester} \\
    (-J_+)W_+(B)+W_+(B)(-J_+^*)+\widetilde{B}_{+,+}&=O\tblue{.} \label{eq:W_unstable_jordan_sylvester}
\end{align}
\end{subequations}

For $i,j=m_-+1,\ldots,m_-+m_0$, we denote the $(n_i,n_j)$-entry of $\widetilde{B}_{i,j}$ by $\widetilde{b}_{i,j}$, which is the lower-right entry, and define
\begin{equation}\label{eq:C_definition}
    C_{i,j}\coloneq
    \begin{pmatrix}
        \frac{1}{(n_i-1)!(n_j-1)!(n_i+n_j-1)} & \cdots & \frac{1}{(n_i-1)!0!n_i} \\
        \vdots & \ddots & \vdots \\
        \frac{1}{0!(n_j-1)!n_j} & \cdots & \frac{1}{0!0!1}
    \end{pmatrix}.
\end{equation}
Moreover, define $W_0^{(\ell)}(B)\in\CC^{n_0^{(\ell)}\times n_0^{(\ell)}} \ (\ell=1,\ldots,m_0')$ and $W_0(B)\in\CC^{n_0\times n_0}$ by
\begin{gather}
    W_0^{(\ell)}(B)\coloneq
    \begin{pmatrix}
        \widetilde{b}_{\alpha+1,\alpha+1}C_{\alpha+1,\alpha+1} & \cdots & \widetilde{b}_{\alpha+1,\beta}C_{\alpha+1,\beta} \\
        \vdots & \ddots & \vdots \\
        \widetilde{b}_{\beta,\alpha+1}C_{\beta,\alpha+1} & \cdots & \widetilde{b}_{\beta,\beta}C_{\beta,\beta}
    \end{pmatrix}, \label{eq:W0l_block} \\
    W_0(B)\coloneq
     \mathrm{diag}(W_0^{(1)}(B),\ldots,W_0^{(m_0')}(B)), \label{eq:W0_block_diagonal}
\end{gather}
where we set $\alpha=m_-+m_0^{(1)}+\cdots+m_0^{(\ell-1)},\ \beta=\alpha+m_0^{(\ell)}$ to simplify the indexing.

Then, the limit of $W(B;T)$ can be expressed in terms of $W_-(B), W_0(B)$, and $W_+(B)$.

\begin{theorem}\label{prop:scaledW_convergence}
    The scaled controllability Gramian $W(B;T)$ converges as $T\to\infty$, and its limit is given by
    \begin{equation}\label{eq:W_block_diagonal}
        W_\infty(B)\coloneq
        \mathrm{diag}(W_-(B),W_0(B),W_+(B)). 
    \end{equation}
\end{theorem}

\begin{proof}
    See Appendix~\ref{sec:proof_scaledW_convergence}.
    \qed
\end{proof}

A closely related polynomial scaling for the controllability
Gramian of a controllable pair with a purely imaginary spectrum
was established in \cite[Theorem~4.1]{belabbas2025infinite}.
Theorem~\ref{prop:scaledW_convergence} differs in that it treats stable, imaginary-axis, and
unstable Jordan blocks simultaneously, allows a general
matrix weight $B$, and identifies the explicit blockwise limit
needed for the subsequent optimization analysis.

Analogously to $W(p;T)$, we abbreviate $W_\infty(\diag(p))$ as $W_\infty(p)$ for $p\in\RR^n$, and refer to it as the infinite-horizon scaled controllability Gramian.
The same notation is used for $W_-, W_0$, and $W_+$.

For the finite-horizon controllability scoring problem, the following lemma guarantees that the feasible set is nonempty.

\begin{lemma}\label{prop:finite_nonempty_lemma}
    Assume that $T>0$ is finite, and let $p\in\Delta$ and $p_i>0 \ (i=1,\ldots,n)$.
    Then, $\widetilde{W}(p;T)\succ O$; equivalently, $W(p;T)\succ O$.
\end{lemma}

In the subsequent sections, we study controllability scoring problems based on matrices such as $W_\infty(p)$.
The following lemma, analogous to \cref{prop:finite_nonempty_lemma}, guarantees that the feasible set is nonempty.

\begin{lemma}\label{prop:infinite_nonempty_lemma}
    Let $p\in\Delta$ and $p_i>0 \ (i=1,\ldots,n)$.
    Then, $W_\infty(p)\succ O$.
    In particular, $W_-(p)\succ O$.
\end{lemma}

\begin{proof}[\cref{prop:finite_nonempty_lemma,prop:infinite_nonempty_lemma}]
    See Appendix~\ref{sec:proof_nonempty_lemma}.
    \qed
\end{proof}

\begin{remark}
Suppose that $A$ is diagonalizable over $\CC$, so that all Jordan
blocks have size one.
Then,
\[
D_i(T)=
\begin{cases}
1& {\rm if}\quad  \re(\lambda_i)<0,\\
T^{1/2}& {\rm if}\quad \re(\lambda_i)=0,\\
\e^{\lambda_iT}& {\rm if}\quad\re(\lambda_i)>0.
\end{cases}
\]
Let $\widetilde B=Q^{-1}BQ^{-*}=(\widetilde b_{ij})$.
For the stable and unstable parts, the entries of the limiting
matrices are
\[
(W_-(B))_{ij}
=
-\frac{\widetilde b_{ij}}
{\lambda_i+\overline{\lambda_j}},
\qquad
(W_+(B))_{ij}
=
\frac{\widetilde b_{ij}}
{\lambda_i+\overline{\lambda_j}},
\]
for the corresponding stable and unstable indices, respectively.
For imaginary-axis eigenvalues,
\[
\frac{1}{T}
\int_0^T
\e^{(\lambda_i+\overline{\lambda_j})t}
\widetilde b_{ij}\,\rd t
\rightarrow
\begin{cases}
\widetilde b_{ij}& {\rm if}\quad\lambda_i=\lambda_j,\\
0& {\rm if}\quad\lambda_i\neq\lambda_j.
\end{cases}
\]
Hence, the center block of the limiting Gramian consists of the
submatrices of $\widetilde B$ associated with equal imaginary-axis
eigenvalues. The more general construction in
\cref{eq:C_definition,eq:W0l_block} accounts for the additional polynomial
growth caused by nontrivial Jordan blocks.
\end{remark}

\subsection{Reformulation of controllability scoring problems}
\label{subsec:reformulation}

We verify that \cref{prob:scaled_vcs,prob:scaled_aecs} are equivalent to \cref{prob:vcs,prob:aecs}, respectively.
It follows from \cref{eq:scaledW_definition} that
\begin{equation*}
    \widetilde{W}(p;T)=QD(T)W(p;T)D(T)^*Q^*
\end{equation*}
for $p\in\RR^n$.
Consequently,
\begin{multline*}
    -\log\det \widetilde{W}(p;T) \\
    =-\log\det W(p;T)-\log\det \left(QD(T)D(T)^*Q^*\right),
\end{multline*}
where the second term is independent of $p\in\RR^n$.
Moreover, since $Q$ and $D(T)$ are invertible, $\widetilde{W}(p;T)\succ O$ and $W(p;T)\succ O$ are equivalent.
Therefore, \cref{prob:scaled_vcs} is equivalent to \cref{prob:vcs}.
Similarly, since
\begin{multline*}
    \tr\left(\widetilde{W}(p;T)^{-1}\right) \\
    =\tr\left(Q^{-*}D(T)^{-*}W(p;T)^{-1}D(T)^{-1}Q^{-1}\right),
\end{multline*}
\cref{prob:scaled_aecs} is equivalent to \cref{prob:aecs}.

Since the uniqueness of the optimal solutions is equivalent as well, we obtain the same result as in \cite[Theorem 5]{SatoTerasaki2024} and \cite[Theorem 1]{SatoKawamura2025}.
In Appendix~\ref{sec:finite_uniqueness}, we summarize existing results on the uniqueness of the optimal solutions and present new results.

If $A$ is not Hurwitz, $\widetilde{W}(p;T)$ diverges, and hence evaluating the objective functions of \cref{prob:vcs,prob:aecs} becomes numerically unstable.
On the other hand, from \cref{prop:scaledW_convergence}, $W(p;T)$ converges to $W_\infty(p)$, and it is positive definite if $p_i>0 \ (i=1,\ldots,n)$.
Therefore, the objective function of \cref{prob:scaled_vcs} remains numerically stable even for large $T$.
Moreover, since $D(T)^{-1}$ converges, the objective function of \cref{prob:scaled_aecs} is also numerically stable.
Hence, challenge 2) in \cref{subsec:challenges} is thereby resolved.

\section{Extension to infinite-horizon settings}
\label{sec:extension}

In this section, we formalize the infinite-horizon controllability scoring problems by taking the limit $T\to\infty$ in \cref{prob:scaled_vcs,prob:scaled_aecs}.
We then justify this formulation by proving that, under suitable assumptions, the optimal solutions are uniquely determined and that the finite-horizon controllability scores indeed converge to them as $T\to\infty$.

In what follows, we define
    $W_{\infty,i}\coloneq W_\infty(e_ie_i^\top)$,
for which
\begin{equation*}
    W_\infty(p)=\sum_{i=1}^n p_i W_{\infty,i}
\end{equation*}
holds for $p=(p_1,\ldots,p_n)^\top\in\RR^n$.
We define $W_{-,i}, W_{0,i}$, and $W_{+,i}$ analogously.

\subsection{Problem formulation}
\label{subsec:infinite_formulation}

\begin{table*}[t]
    \centering
    \caption{Scope of the limiting formulations and their real-Schur implementations.}
    \label{tab:scope_summary}
    \begin{tabular}{p{0.18\textwidth}p{0.30\textwidth}p{0.43\textwidth}}
        \hline
        Formulation
        & Spectral requirement
        & Spectral information retained in the objective
        \\ \hline\hline
        Infinite-horizon VCS (Problem~\ref{prob:infinite_vcs})
        & Arbitrary $A$; neither Assumption~\ref{assump:stable_eigenvalue}
        nor Assumption~\ref{assump:imaginary_eigenvalue}
        & Stable, center, and unstable components
        \\
        Infinite-horizon AECS (Problem~\ref{prob:infinite_aecs})
        & At least one stable eigenvalue
        (Assumption~\ref{assump:stable_eigenvalue});
        Assumption~\ref{assump:imaginary_eigenvalue} is not required
        & Stable component only
        \\
        Real-Schur Problems~\ref{prob:schur_vcs}
        and~\ref{prob:schur_aecs}
        & Assumption~\ref{assump:imaginary_eigenvalue},
        and additionally Assumption~\ref{assump:stable_eigenvalue}
        for AECS
        & The same components as the corresponding limiting problem;
        the restriction concerns computation, not the limiting theory
        \\ \hline
    \end{tabular}
\end{table*}

Table~\ref{tab:scope_summary} summarizes the spectral requirements
of the two limiting problems and distinguishes them from the
additional requirement imposed only for the real-Schur
implementation in \cref{sec:algorithms}.
Both limiting formulations allow arbitrary imaginary-axis eigenvalues
and associated Jordan blocks.
The infinite-horizon VCS applies to an arbitrary $A$ and retains all
spectral components, whereas the infinite-horizon AECS requires a
nonempty stable component and retains only that component.
Assumption~\ref{assump:imaginary_eigenvalue} in Section~\ref{sec:algorithms} therefore concerns
computation rather than the limiting theory.

We first formulate the infinite-horizon version of \cref{prob:scaled_vcs} as follows:
\begin{problem}\label{prob:infinite_vcs}
    \begin{align*}
        \minimize_{p} \quad & {-}\log\det W_\infty(p) \\
        \mathrm{subject\ to} \quad & p\in\Delta,\ W_\infty(p)\succ O.
    \end{align*}
\end{problem}
We refer to the optimal solution to \cref{prob:infinite_vcs} as the VCS on an infinite time horizon, or simply the infinite-horizon VCS.
Since $W(p;T)\to W_\infty(p)$ as $T\to\infty$, \cref{prob:infinite_vcs} is obtained from \cref{prob:scaled_vcs} by taking the limit $T\to\infty$, and its objective depends on the full matrix $W_\infty(p)$.

To formulate the infinite-horizon AECS, we require the stable block
to be nonempty.
We therefore impose the following assumption.

\begin{assumption}\label{assump:stable_eigenvalue}
    The matrix $A$ has at least one eigenvalue with a negative real part, i.e., $n_-\geq1$.
\end{assumption}

Under \cref{assump:stable_eigenvalue}, we formulate the infinite-horizon version of \cref{prob:scaled_aecs} as follows:

\begin{problem}\label{prob:infinite_aecs}
    \begin{align*}
        \minimize_{p} \quad & \tr\left(Q^{-*}
         \mathrm{diag}(W_-(p)^{-1},O,O)
        Q^{-1}\right) \\
        \mathrm{subject\ to} \quad & p\in\Delta,\ W_-(p)\succ O.
    \end{align*}
\end{problem}

For every $p$ satisfying $W_\infty(p)\succ O$,
\cref{prop:scaledW_convergence} and continuity of matrix inversion give
    $W(p;T)^{-1}\to W_\infty(p)^{-1}$.
Moreover,
    $D(T)^{-1}\to \diag(I_{n_-},O,O)$.
Since
    $W_\infty(p)
    =
    \diag\left(W_-(p),W_0(p),W_+(p)\right)$,
we obtain
\begin{equation}
\label{eq:aecs_pointwise_matrix_limit}
    D(T)^{-*}W(p;T)^{-1}D(T)^{-1}
    \to
    \diag\left(W_-(p)^{-1},O,O\right).
\end{equation}
Hence, only the stable block contributes to the limiting AECS
objective, which explains the form of
\cref{prob:infinite_aecs}.
The convergence of the optimization problems and their optimal
solutions is established rigorously through epi-convergence in
\cref{prop:aecs_convergence}.

We refer to an optimal solution of
\cref{prob:infinite_aecs} as an AECS on an infinite time
horizon, or simply an infinite-horizon AECS.
If $A$ is Hurwitz, then \cref{prob:infinite_vcs,prob:infinite_aecs} coincide with the controllability scoring problems proposed in \cite{SatoTerasaki2024}.

\begin{remark}\label{rem:vcs_aecs_controllability}
    \cref{prob:infinite_vcs} imposes controllability of the entire
    system, whereas \cref{prob:infinite_aecs} requires controllability
    only of the stable modes.
    Consequently, the infinite-horizon VCS preserves full-system
    controllability, whereas the infinite-horizon AECS may assign zero
    weight to an input that is indispensable for a center or unstable
    mode.
    This spectral distinction is consistent with the finite-horizon structural analysis in \cite[Section~IV]{SatoOED2026}, where a lower bound on the AECS objective is derived to characterize the long-horizon downweighting of certain source-like nodes.
    The relation between that analysis and the zero AECS allocation
    observed in the directed-Laplacian example is discussed in
    Section~\ref{sec:experiments}.
\end{remark}

\subsection{Convergence of controllability scores}
\label{subsec:convergence}

\Cref{prob:infinite_vcs,prob:infinite_aecs} are obtained by formally taking the limit $T\to\infty$ in \cref{prob:vcs,prob:aecs}.
Under suitable assumptions, the optimal solutions to \cref{prob:vcs,prob:aecs} (finite-horizon controllability scores) indeed converge to those to \cref{prob:infinite_vcs,prob:infinite_aecs} (infinite-horizon controllability scores). 

\begin{theorem}\label{prop:vcs_convergence}
     Assume that \cref{prob:infinite_vcs} admits a unique optimal solution
    $p_{\mathrm{VCS}}^{\infty}$.
    For each $T>0$, let $p_{\mathrm{VCS}}^{T}$ be any optimal solution to
    \cref{prob:vcs}.
    Then, $p_{\mathrm{VCS}}^{T}
        \to
        p_{\mathrm{VCS}}^{\infty}$ as $T\to\infty$.
\end{theorem}

\begin{theorem}\label{prop:aecs_convergence}
      Assume that \cref{assump:stable_eigenvalue} holds and that
    \cref{prob:infinite_aecs} admits a unique optimal solution
    $p_{\mathrm{AECS}}^{\infty}$.
    For each $T>0$, let $p_{\mathrm{AECS}}^{T}$ be any optimal solution to
    \cref{prob:aecs}.
    Then, $p_{\mathrm{AECS}}^{T}
        \to
        p_{\mathrm{AECS}}^{\infty}$ as $T\to\infty$.
\end{theorem}

\begin{proof}[\cref{prop:vcs_convergence,prop:aecs_convergence}]
    See Appendix~\ref{sec:convergence_proof}.
    \qed
\end{proof}

Theorems~\ref{prop:vcs_convergence} and
\ref{prop:aecs_convergence} use uniqueness only to obtain convergence
to a single limiting allocation.
Without uniqueness, every cluster point is still an infinite-horizon
optimizer, and the distance to the limiting optimal-solution set
converges to zero, as stated next.

\begin{proposition}
\label{cor:nonunique_convergence}
    Let $\mathcal P_{\mathrm{VCS}}^\infty$ denote the set of optimal solutions to \cref{prob:infinite_vcs}, and for every $T>0$ select any optimal solution $p_{\mathrm{VCS}}^T$ to \cref{prob:vcs}.  Every cluster point of $p_{\mathrm{VCS}}^T$ as $T\to\infty$ belongs to $\mathcal P_{\mathrm{VCS}}^\infty$, and
    \begin{equation*}
        \inf_{q\in\mathcal P_{\mathrm{VCS}}^\infty}\norm{p_{\mathrm{VCS}}^T-q}\rightarrow0.
    \end{equation*}
    Under \cref{assump:stable_eigenvalue}, the analogous statements hold for the AECS and its optimal-solution set $\mathcal P_{\mathrm{AECS}}^\infty$.
\end{proposition}

\begin{proof}
    See Appendix~\ref{sec:convergence_proof}.
    \qed
\end{proof}

\Cref{prop:vcs_convergence,prop:aecs_convergence} assume the uniqueness of the optimal solutions to \cref{prob:infinite_vcs,prob:infinite_aecs}.
As discussed in \cref{subsec:infinite_uniqueness}, this assumption is satisfied under appropriate conditions; otherwise, \cref{cor:nonunique_convergence} gives convergence to the solution set but the selected minimizers need not converge to one point.

\subsection{Uniqueness analysis}
\label{subsec:infinite_uniqueness}

We discuss uniqueness of the optimal solutions to \cref{prob:infinite_vcs,prob:infinite_aecs}, which is important for two reasons.

First, the convergence results for the finite-horizon controllability scores established in \cref{subsec:convergence} are derived under the assumption that the optimal solutions to \cref{prob:infinite_vcs,prob:infinite_aecs} are unique.
When the optimal solutions are unique, the finite-horizon controllability scores converge to them as $T\to\infty$.

Second, uniqueness is essential from the viewpoints of reproducibility and interpretability.
As has been studied in detail for finite-horizon controllability scores~\cite{SatoTerasaki2024,SatoKawamura2025}, the uniqueness is crucial when they are used as centrality measures.
Existing results on the uniqueness of finite-horizon controllability scores are summarized in Appendix \ref{sec:finite_uniqueness}.
If the optimal solution is not unique, different researchers may obtain different infinite-horizon controllability scores and consequently draw different conclusions, which undermines reproducibility.
Moreover, the presence of multiple optimal solutions raises interpretability issues regarding how such solutions should be used to define a meaningful centrality measure.

We show that although the uniqueness of the optimal solution depends on the structure of the eigenvectors of the matrix $A$, it holds in a wide range of practically relevant cases.
Note that even if the finite-horizon controllability score is unique for all $T>0$, the infinite-horizon controllability score need not be unique.
Indeed, the following example shows that the finite-horizon controllability scores are unique for all $T>0$, whereas the infinite-horizon controllability scores are not unique.

\begin{example}\label{ex:infinite_nonuniqueness}
    Let $A=\begin{pmatrix} 0 & 1 \\ 1 & 0 \end{pmatrix}$.
    Then, \cref{prob:vcs,prob:aecs} each admit a unique optimal solution for all $T>0$.
    However, \cref{prob:infinite_vcs,prob:infinite_aecs} have multiple optimal solutions.

    Indeed, by taking $Q=\dfrac{1}{\sqrt{2}}\begin{pmatrix} 1 & -1 \\ 1 & 1 \end{pmatrix}$, we obtain $Q^{-1}AQ=\diag(1,-1)$, whose eigenvalues are $1$ and $-1$.
    Therefore, $\Theta'=\emptyset$ in \cref{prop:uniqueness_by_eigenvalue} of Appendix~\ref{sec:finite_uniqueness}, and hence, by \cref{prop:uniqueness_by_eigenvalue}, \cref{prob:vcs,prob:aecs} each admit a unique optimal solution for all $T>0$.
    
    On the other hand, the infinite-horizon scaled controllability Gramians, are given by $W_{\infty,1}=W_{\infty,2}=\begin{pmatrix} \frac{1}{4} & 0 \\ 0 & \frac{1}{4} \end{pmatrix}$ and $W_{-,1}=W_{-,2}=\frac{1}{4}$.
    Thus, $W_\infty(p)$ and $W_-(p)$ take constant matrices for all $p\in\Delta$.
    As a consequence, every point in $\Delta$ is an optimal solution to \cref{prob:infinite_vcs,prob:infinite_aecs}. \qed
\end{example}

We partition $Q^{-1}e_i$ into blocks in following two ways:
    $Q^{-1}e_i=
    \begin{pmatrix}
        q_{1,i} \\
        \vdots \\
        q_{m,i}
    \end{pmatrix}
    =
    \begin{pmatrix}
        q_{-,i} \\
        q_{0,i} \\
        q_{+,i}
    \end{pmatrix}$.
Furthermore, let $q'_{k,i}$ denote the $n_k$-th (i.e., the last) entry of $q_{k,i}$.
For $\ell=1,\ldots,m_0'$, we define
    $q_{0,i}^{(\ell)}\coloneq(q'_{\alpha+1,i},\ldots,q'_{\beta,i})^\top\in\CC^{m_0^{(\ell)}}$,
where we set $\alpha=m_-+m_0^{(1)}+\cdots+m_0^{(\ell-1)},\ \beta=\alpha+m_0^{(\ell)}$ to simplify the indexing.
Then, the following theorems hold.

\begin{theorem}\label{prop:infinite_vcs_uniqueness}
    Let
    \begin{align*}
        M_{\mathrm{VCS}}&\coloneq
        \begin{pmatrix}
            q_{-,1}\otimes q_{-,1} & \cdots & q_{-,n}\otimes q_{-,n} \\
            q_{0,1}^{(1)}\otimes q_{0,1}^{(1)} & \cdots & q_{0,n}^{(1)}\otimes q_{0,n}^{(1)} \\
            \vdots & \ddots & \vdots \\
            q_{0,1}^{(m_0')}\otimes q_{0,1}^{(m_0')} & \cdots & q_{0,n}^{(m_0')}\otimes q_{0,n}^{(m_0')} \\
            q_{+,1}\otimes q_{+,1} & \cdots & q_{+,n}\otimes q_{+,n}
        \end{pmatrix} \\
        &\in\CC^{(n_-^2+(m_0^{(1)})^2+\cdots+(m_0^{(m_0')})^2+n_+^2)\times n},
    \end{align*}
    and assume $\rank M_{\mathrm{VCS}}=n$.
    Then, \cref{prob:infinite_vcs} admits a unique optimal solution.
\end{theorem}

\begin{theorem}\label{prop:infinite_aecs_uniqueness}
    Let
    \begin{equation*}
        M_{\mathrm{AECS}}\coloneq
        \begin{pmatrix}
            q_{-,1}\otimes q_{-,1} & \cdots & q_{-,n}\otimes q_{-,n}
        \end{pmatrix}
        \in\CC^{n_-^2\times n},
    \end{equation*}
    and assume $\rank M_{\mathrm{AECS}}=n$.
    Then, \cref{prob:infinite_aecs} admits a unique optimal solution.
\end{theorem}

\begin{proof}[\cref{prop:infinite_vcs_uniqueness,prop:infinite_aecs_uniqueness}]
    See Appendix~\ref{sec:infinite_uniqueness_proof}.
    \qed
\end{proof}

The rank tests also identify mechanisms by which the sufficient
uniqueness conditions can fail.
A necessary dimensional requirement for the AECS condition is
$n_-^2\geq n$, whereas for the VCS condition it is
$n_-^2+\sum_{\ell=1}^{m_0'}(m_0^{(\ell)})^2+n_+^2\geq n$.
Even when the row dimension is adequate, symmetry or dynamically
indistinguishable candidate inputs can create linearly dependent
columns, as in \cref{ex:infinite_nonuniqueness}.
Because these conditions are sufficient but not necessary, their
failure does not imply nonuniqueness of the optimizer.

\section{Algorithms}
\label{sec:algorithms}

This section provides a computational realization of the limiting
problems rather than a new optimization algorithm.
As summarized in Table~\ref{tab:scope_summary} in Section~\ref{sec:extension}, the theory in
\cref{sec:reformulation,sec:extension} allows arbitrary
imaginary-axis eigenvalues and associated Jordan blocks.
Once the limiting candidate matrices are available, the resulting
open-domain convex problems can, in principle, be solved by the
Armijo projected-gradient framework in
\cite{SatoOpenDomain2026}.
The remaining numerical difficulty for a general center spectrum is
the stable construction of these matrices without Jordan-chain
information.

The Jordan form is used in the preceding sections only as an
analytical device.
Because its direct numerical recovery is intrinsically ill
conditioned \cite{GolubWilkinson1976}, we instead use real-Schur,
Sylvester, and Lyapunov computations under the following additional
assumption.

\begin{assumption}\label{assump:imaginary_eigenvalue}
 The only possible eigenvalue of $A$ on the imaginary axis is zero,
    and zero, if present, is semisimple; that is, $A$ has no nonzero
    purely imaginary eigenvalues, and all Jordan blocks associated
    with zero are of size one.
\end{assumption}

Assumption~\ref{assump:imaginary_eigenvalue} restricts only the
matrix construction in this section.
It is used in \cref{prop:infinite_equivalence} and the subsequent
computational procedure, but not in the limiting formulations,
convergence results, or uniqueness analysis.
It includes every matrix with no imaginary-axis eigenvalues, as well
as directed graph Laplacians.
For a more general center spectrum, the limiting theory remains valid,
and the optimization framework remains applicable once the limiting
candidate matrices are available.
A numerically reliable construction of these matrices without
Jordan-chain information, however, is left open.

An important example satisfying this assumption is the graph Laplacian.
\begin{example}\label{prop:laplacian_eigenvalue}
    The matrix $L=(\ell_{ij})\in\RR^{n\times n}$ is said to be the graph Laplacian if, for all $i\neq j$,
    \begin{equation*}
        \ell_{ij}=
        \begin{cases}
            -c_{ij} &
            \begin{array}{l}
                \text{if there is an edge of weight $c_{ij}>0$} \\
                \text{from node $x_j$ to node $x_i$}
            \end{array}, \\
            0 & \text{otherwise},
        \end{cases}
    \end{equation*}
    and the diagonal entries satisfy
        $\ell_{ii}=-\sum_{j\neq i} \ell_{ij}$.
    When considering graph Laplacians, we restrict attention to positive edge weights whenever an edge exists.
    For undirected graphs, the Laplacian matrix $L$ is symmetric, whereas for directed graphs, $L$ is generally not symmetric.

    Let $L$ be the Laplacian matrix of a directed graph.
    Then, $A=-L$ satisfies \cref{assump:imaginary_eigenvalue}.
    In fact, from the definition, the matrix $A=(a_{ij})$ satisfies
       $a_{ij}\geq 0 \ (i\neq j)$ and $a_{ii}=-\sum_{j\neq i}a_{ij}\leq 0$.
    It follows from Ger\v shgorin's theorem \cite[Theorem 1.1]{Varga2004} that all eigenvalues of $A$ have nonpositive real parts, and that zero is the only eigenvalue that may lie on the imaginary axis.
    Furthermore, $\rank(A)=\rank(A^2)$ holds~\cite[Proposition 12]{ChebotarevAgaev2002}.
    This implies that the zero eigenvalue of $A$ is semisimple. \qed
\end{example}

First, we perform a block diagonalization using the real Schur decomposition and the solution of Sylvester equations.
This approach avoids the intrinsically ill-conditioned numerical recovery of the Jordan canonical form.
By the real Schur decomposition, we can find an orthogonal matrix $R^{(1)}\in\RR^{n\times n}$ such that
\begin{equation}\label{eq:A_schur_decomposition}
    \left(R^{(1)}\right)^\top AR^{(1)}=
    \begin{pmatrix}
        A_- & O & O \\
        A_{0,-} & A_0 & O \\
        A_{+,-} & A_{+,0} & A_+
    \end{pmatrix}\tblue{,}
\end{equation}
where $A_-,A_0$, and $A_+$ have eigenvalues with negative, zero, and positive real parts, respectively.
While the real Schur decomposition typically yields a block upper-triangular matrix, applying it to $A^\top$ yields a block lower-triangular matrix.

Let $R^{(2)}_{0,-}, R^{(2)}_{+,-}$, and $R^{(2)}_{+,0}$ be matrices satisfying
\begin{subequations}\label{eq:schur_sylvester_group}
\begin{align}
    A_0 R^{(2)}_{0,-}-R^{(2)}_{0,-}A_-&=-A_{0,-}, \label{eq:schur_sylvester_1}\\
    A_+ R^{(2)}_{+,0}-R^{(2)}_{+,0}A_0&=-A_{+,0}, \label{eq:schur_sylvester_2}\\
    A_+ R^{(2)}_{+,-}-R^{(2)}_{+,-}A_-&=-(A_{+,-}+A_{+,0}R^{(2)}_{0,-})\tblue{.} \label{eq:schur_sylvester_3}
\end{align}
\end{subequations}
Since no two of the matrices $A_-,A_0$, and $A_+$ share a common eigenvalue, \cref{eq:schur_sylvester_1,eq:schur_sylvester_2,eq:schur_sylvester_3} can be solved successively and admit a unique solution from \cref{prop:sylvester_uniqueness}.
Define
    $R^{(2)}\coloneq
    \begin{pmatrix}
        I_{n_-} & O & O \\
        R^{(2)}_{0,-} & I_{n_0} & O \\
        R^{(2)}_{+,-} & R^{(2)}_{+,0} & I_{n_+}
    \end{pmatrix}$
and set $R=R^{(1)}R^{(2)}$.
Then, by using
\begin{equation}\label{eq:S_inverse}
    \left(R^{(2)}\right)^{-1}=
    \begin{pmatrix}
        I_{n_-} & O & O \\
        -R^{(2)}_{0,-} & I_{n_0} & O \\
        -R^{(2)}_{+,-}+R^{(2)}_{+,0}R^{(2)}_{0,-} & -R^{(2)}_{+,0} & I_{n_+}
    \end{pmatrix}\tblue{,}
\end{equation}
we obtain the block-diagonal form
\begin{equation}\label{eq:A_block_diagonal}
    R^{-1}AR= \mathrm{diag} (A_-, A_0, A_+).
\end{equation}

Then, let
\begin{equation}\label{eq:q_definition}
    \begin{pmatrix}
        r_{-,i} \\
        r_{0,i} \\
        r_{+,i}
    \end{pmatrix}
    \coloneq R^{-1}e_i,
\end{equation}
where $r_{-,i}\in\RR^{n_-},\ r_{0,i}\in\RR^{n_0},\ r_{+,i}\in\RR^{n_+}$, and let $V_{-,i}\in\RR^{n_-\times n_-}$ and $V_{+,i}\in\RR^{n_+\times n_+}$ be the solutions to
\begin{subequations}\label{eq:W_schur_sylvester}
\begin{align}
    A_- V_{-,i}+V_{-,i}A_-^\top+r_{-,i}r_{-,i}^\top&=O, \label{eq:W_stable_schur_sylvester}\\
    (-A_+) V_{+,i}+V_{+,i}(-A_+^\top)+r_{+,i}r_{+,i}^\top&=O. \label{eq:W_unstable_schur_sylvester}
\end{align}
\end{subequations}
Moreover, let
\begin{equation}
    V_{0,i}\coloneq r_{0,i}r_{0,i}^\top. \label{eq:W_imaginary_schur}
\end{equation}
Finally, for $p=(p_1,\ldots,p_n)^\top\in\RR^n$, we define
\begin{gather*}
    V_\alpha(p)\coloneq\sum_{i=1}^n p_i V_{\alpha,i} \quad (\alpha\in\{-,0,+\}), \\
     V_i\coloneq \mathrm{diag} (V_{-,i}, V_{0,i}, V_{+,i}), \\
    V(p)\coloneq
    \mathrm{diag} (V_{-}(p), V_{0}(p), V_{+}(p)).
\end{gather*}

We then consider the following alternative problems.
\begin{problem}\label{prob:schur_vcs}
    \begin{align*}
        \minimize_{p} \quad & {-}\log\det V(p) \\
        \mathrm{subject\ to} \quad & p\in\Delta,\ V(p)\succ O.
    \end{align*}
\end{problem}
\begin{problem}\label{prob:schur_aecs}
    \begin{align*}
        \minimize_{p} \quad & \tr\left(V_-(p)^{-1}\right) \\
        \mathrm{subject\ to} \quad & p\in\Delta,\ V_-(p)\succ O.
    \end{align*}
\end{problem}

In fact, \cref{prob:infinite_vcs,prob:infinite_aecs} are equivalent to \cref{prob:schur_vcs,prob:schur_aecs}, respectively.

\begin{theorem}\label{prop:infinite_equivalence}
    Under \cref{assump:imaginary_eigenvalue}, \cref{prob:infinite_vcs,prob:schur_vcs} are equivalent in the sense that they share the same set of optimal solutions.
    Under \cref{assump:stable_eigenvalue,assump:imaginary_eigenvalue}, \cref{prob:infinite_aecs,prob:schur_aecs} are also equivalent.
\end{theorem}

\begin{proof}
    See Appendix~\ref{sec:infinite_equivalence_proof}.
    \qed
\end{proof}

Although \cref{prob:infinite_vcs,prob:infinite_aecs} are formulated based on the Jordan canonical form, under \cref{assump:stable_eigenvalue,assump:imaginary_eigenvalue}, they are respectively equivalent to \cref{prob:schur_vcs,prob:schur_aecs}.
Hence, it suffices to compute the optimal solutions to \cref{prob:schur_vcs,prob:schur_aecs}, without resorting to the Jordan canonical form.

\begin{algorithm}[t]
    \begin{algorithmic}[1]
        \REQUIRE The matrix $A$.
        \STATE Perform the real Schur decomposition~\cref{eq:A_schur_decomposition} and compute $R^{(1)}$.
        \STATE Compute $R^{(2)}_{0,-}, R^{(2)}_{+,-}$, and $R^{(2)}_{+,0}$ by \cref{eq:schur_sylvester_1,eq:schur_sylvester_2,eq:schur_sylvester_3}, and then compute $(R^{(2)})^{-1}$ by \cref{eq:S_inverse}.
        \STATE Compute $r_{-,i}, r_{0,i}$, and $r_{+,i}$ by \cref{eq:q_definition}.
        \STATE Compute $V_{-,i}$ by solving \cref{eq:W_stable_schur_sylvester}.
        \STATE $V_{0,i}=r_{0,i}r_{0,i}^\top$.
        \STATE Compute $V_{+,i}$ by solving \cref{eq:W_unstable_schur_sylvester}.
        \ENSURE $V_{-,i}, V_{0,i}$, and $V_{+,i}$
    \end{algorithmic}
    \caption{The computation of $V_{-,i}, V_{0,i}$, and $V_{+,i}$\tblue{.}}
    \label{alg:schur_computation}
\end{algorithm}

To summarize the discussion so far, under \cref{assump:stable_eigenvalue,assump:imaginary_eigenvalue}, $V_{-,i}, V_{0,i}$, and $V_{+,i}$ can be computed by \cref{alg:schur_computation}.
After computing $V_{-,i}$, $V_{0,i}$, and $V_{+,i}$
by \cref{alg:schur_computation}, we solve
\cref{prob:schur_vcs,prob:schur_aecs} using the domain-aware
Armijo projected-gradient method in
\cite{SatoTerasaki2024,SatoOpenDomain2026}.
Let
\begin{equation*}
    g(p)\coloneq-\log\det V(p),
    \qquad
    h(p)\coloneq\tr\left(V_-(p)^{-1}\right).
\end{equation*}
Their gradients are
\begin{align*}
    (\nabla g(p))_i
    &=-\tr\left(V(p)^{-1}V_i\right),\\
    (\nabla h(p))_i
    &=-\tr\left(
        V_-(p)^{-1}V_{-,i}V_-(p)^{-1}
    \right).
\end{align*}
For either objective $f\in\{g,h\}$, the iteration has the form
\begin{equation*}
    p^{(k+1)}
    =
    \Pi_{\Delta}\left(
        p^{(k)}-\alpha^{(k)}\nabla f(p^{(k)})
    \right),
\end{equation*}
where $\alpha^{(k)}$ is selected by domain-aware Armijo
backtracking so that the accepted iterate remains in the
positive-definite domain.
We initialize the method at the uniform allocation
$p^{(0)}=(1/n,\ldots,1/n)^\top$, which treats all candidate
inputs symmetrically and is strictly feasible by
\cref{prop:infinite_nonempty_lemma,prop:infinite_equivalence}.
The initialization may affect the iteration path and convergence
speed, and, when the optimizer is nonunique, which optimizer is
obtained, but it does not change the optimal set.
Well-definedness and convergence follow from the cited results.

Following \cite[Section~IV]{umezu2026controllability}, the dense implementation
requires $O(n^3)$ storage and $O(n^3)$ work per projected-gradient
iteration.
Finite-horizon quadrature costs $O(N_Tn^3)$, whereas direct dense
construction of all candidate matrices, including the present
real-Schur/Lyapunov procedure, has $O(n^4)$ worst-case cost but is
independent of $T$.
Thus, neither formulation is uniformly cheaper.
For large sparse systems, more recent low-rank ADI variants with
improved shift selection and rational Krylov subspace methods may
reduce both preprocessing cost and storage
\cite{Simoncini2007,BennerKurschnerSaak2014}.
An optimization method operating directly on the resulting low-rank
factors is left for future work.

\section{Numerical experiments}
\label{sec:experiments}

\begin{figure}[htbp]
    \centering
    \begin{tikzpicture}[scale=0.7]
        \node (n1) at (4.1,2.3) {1};
        \node (n2) at (0.5,2.3) {2};
        \node (n3) at (2.9,2.3) {3};
        \node (n4) at (1.7,2.3) {4};
        \node (n5) at (4.1,0.0) {5};
        \node (n6) at (1.7,0) {6};
        \node (n7) at (2.3,3.5) {7};
        \node (n8) at (2.9,0.0) {8};
        \node (n9) at (5.1,3.5) {9};
        \node (n10) at (0.5,1.0) {10};
        
        \draw [thick] (n1) circle [radius=0.28];
        \draw [thick] (n2) circle [radius=0.28];
        \draw [thick] (n3) circle [radius=0.28];
        \draw [thick] (n4) circle [radius=0.28];
        \draw [thick] (n5) circle [radius=0.28];
        \draw [thick] (n6) circle [radius=0.28];
        \draw [thick] (n7) circle [radius=0.28];
        \draw [thick] (n8) circle [radius=0.28];
        \draw [thick] (n9) circle [radius=0.28];
        \draw [thick] (n10) circle [radius=0.28];

        \draw [-Latex,thick] (n1)--(n5);
        \draw [-Latex,thick] (n2)--(n10);
        \draw [-Latex,thick] (n3)--(n8);
        \draw [-Latex,thick] (n4)--(n6);
        \draw [-Latex,thick] (n7)--(n1);
        \draw [-Latex,thick] (n7)--(n2);
        \draw [-Latex,thick] (n7)--(n3);
        \draw [-Latex,thick] (n7)--(n4);
        \draw [-Latex,thick] (n9)--(n1);
        \draw [-Latex,thick] (n10)--(n6);
    \end{tikzpicture}
    \caption{Directed network employed in the numerical experiments.}
    \label{fig:network}
\end{figure}
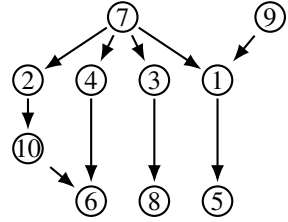

In this section, we conduct numerical experiments using the same network as in \cite{SatoKawamura2025}, which is illustrated in \cref{fig:network}.
The network consists of $n=10$ nodes, and all edges have a uniform weight $0.2$.
Let $L$ denote the graph Laplacian of this network, and set $A=-L$.

While \cite{SatoKawamura2025} reports finite-horizon
controllability scores up to $T=10000$, we compute the corresponding
infinite-horizon scores using
\cref{prob:schur_vcs,prob:schur_aecs}.
Table~\ref{tab:horizon_comparison} compares the allocations at
$T=10000$ and $T=\infty$.
The maximum componentwise differences are
$1.5\times10^{-4}$ for the VCS and
$2.34\times10^{-3}$ for the AECS, which is consistent with the
convergence established in
\cref{prop:vcs_convergence,prop:aecs_convergence}.
The largest remaining AECS difference occurs at node~9, whose
allocation approaches zero as $T$ increases.
For this example, the matrices appearing in \cref{prop:infinite_vcs_uniqueness,prop:infinite_aecs_uniqueness} satisfied
    $\rank M_{\mathrm{VCS}}
    =
    \rank M_{\mathrm{AECS}}
    =
    10$.
Hence,
\cref{prop:infinite_vcs_uniqueness,prop:infinite_aecs_uniqueness}
certify that both infinite-horizon allocations are unique.

\begin{table}[t]
    \centering
    \caption{Comparison of the controllability scores at
    $T=10000$ and $T=\infty$.}
    \label{tab:horizon_comparison}
    \begin{tabular}{c|cc|cc}
        \hline
        & \multicolumn{2}{c|}{VCS}
        & \multicolumn{2}{c}{AECS} \\
        Node
        & $T=10000$ & $T=\infty$
        & $T=10000$ & $T=\infty$
        \\ \hline\hline
        1  & 0.07327 & 0.07329 & 0.17281 & 0.17428 \\
        2  & 0.10108 & 0.10104 & 0.11364 & 0.11408 \\
        3  & 0.10874 & 0.10871 & 0.12093 & 0.12174 \\
        4  & 0.08636 & 0.08634 & 0.10610 & 0.10610 \\
        5  & 0.04499 & 0.04490 & 0.09230 & 0.09275 \\
        6  & 0.06071 & 0.06069 & 0.13383 & 0.13433 \\
        7  & 0.24952 & 0.24967 & 0.09275 & 0.09233 \\
        8  & 0.04221 & 0.04220 & 0.06945 & 0.06901 \\
        9  & 0.16674 & 0.16686 & 0.00234 & 0.00000 \\
        10 & 0.06632 & 0.06631 & 0.09586 & 0.09538 \\
        \hline
    \end{tabular}
\end{table}

The two limiting allocations differ in whether they preserve
full-system controllability.
To verify this distinction, define
\begin{align*}
    \mathcal C(p)
    &\coloneq
    \begin{bmatrix}
        B(p)&AB(p)&\cdots&A^{n-1}B(p)
    \end{bmatrix},\\
    B(p)
    &\coloneq
    \diag(\sqrt{p_1},\ldots,\sqrt{p_n}).
\end{align*}
For the infinite-horizon VCS allocation,
$\rank\mathcal C(p)=10$, whereas the infinite-horizon AECS assigns
$p_9=0$ and yields $\rank\mathcal C(p)=9$.

To identify the lost controllable direction, let
\begin{equation*}
    B_{-9}
    \coloneq
    \begin{bmatrix}
        e_1&\cdots&e_8&e_{10}
    \end{bmatrix}.
\end{equation*}
The eigenvalues of $A$ are $-0.4$, $-0.2$, and $0$, with
algebraic multiplicities $2$, $6$, and $2$, respectively, and the
zero eigenvalue is semisimple.
A direct PBH test gives
\begin{align*}
    \rank
    \begin{bmatrix}
        \lambda I_n-A&B_{-9}
    \end{bmatrix}
    &=
    \begin{cases}
        10& {\rm if}\quad   \lambda=-0.4,-0.2,\\
        9& {\rm if}\quad   \lambda=0.
    \end{cases}
\end{align*}
Hence, without the input at node~9, all stable modes remain
controllable, whereas exactly one direction in the two-dimensional
center subspace becomes uncontrollable.
The input at node~9 is therefore indispensable for controlling this
center direction.

Node~9 is a source-like node without a negative self-loop, and its
decreasing finite-horizon AECS allocation is consistent with the
long-horizon structural analysis in
\cite[Section~IV]{SatoOED2026}.
Consequently, the VCS retains the contribution of node~9 and assigns
$p_9=0.16686$, whereas the AECS evaluates only the stable block and
assigns $p_9=0$.
This verifies the mechanism described in
\cref{rem:vcs_aecs_controllability}.

\begin{remark}\label{rem:source_code}
    The MATLAB package used to compute the controllability scores is
    publicly available in \cite{ControllabilityScoreCode}.
    The package is not specific to the network considered here: it
    accepts a user-specified square system matrix $A$ and computes
    finite-horizon VCS and AECS, as well as their infinite-horizon
    counterparts within the implementation scope described in
    \cref{sec:algorithms}.
    It can therefore be used to reproduce the reported score
    calculations and to examine other systems.
\end{remark}

\section{Conclusion}
\label{sec:conclusion}

We formulated infinite-horizon VCS and AECS for linear
time-invariant systems whose conventional controllability Gramians may
diverge.
A mode-dependent scaling preserves the finite-horizon optimal
allocations exactly while yielding a convergent scaled Gramian.
Using epi-convergence, we established convergence to the limiting
solution set and derived explicit sufficient rank conditions for
uniqueness.

The limiting VCS and AECS reflect distinct controllability criteria,
consistent with their optimal-experimental-design interpretations in
\cite{SatoOED2026}.
The VCS retains all spectral components and preserves full-system
controllability, whereas the AECS depends only on the stable component
and may neglect an input required for a center or unstable mode.
In the directed-Laplacian example, this distinction resulted in
controllability ranks $10$ and $9$, respectively.

Under an additional condition on the center spectrum, the limiting
problems can be computed using real Schur decompositions and Sylvester
and Lyapunov equations.
Sparse and low-rank implementations and applications to real networks
with non-Hurwitz identified dynamics remain important directions for
future work.

\appendices

\section{Proof of \cref{prop:scaledW_convergence}}
\label{sec:proof_scaledW_convergence}

To prove \cref{prop:scaledW_convergence}, we prepare the following lemma.

\begin{lemma}[{\cite[Theorem 5]{Lancaster1970}}]\label{prop:sylvester_convergence}
    Let $C_1\in\CC^{\ell_1\times\ell_1}, C_2\in\CC^{\ell_2\times\ell_2}$, and $P\in\CC^{\ell_1\times\ell_2}$.
    If $\re(\mu_1+\mu_2)<0$ holds for any eigenvalues $\mu_1$ of $C_1$ and $\mu_2$ of $C_2$, then
        $\int_0^\infty \e^{C_1 t}P\e^{C_2 t}\rd t$
    converges, and it is the unique solution $Y\in\CC^{\ell_1\times\ell_2}$ satisfying
        $C_1Y+YC_2+P=0$.
    \end{lemma}

\begin{proof}[\cref{prop:scaledW_convergence}]
    We prove that $W(B;T)$ converges to $W_\infty(B)$ as $T\to\infty$.
    By definition of $\widetilde{B}$, we obtain
    \begin{equation}\label{eq:W_jordan}
        \widetilde{W}(B;T)=Q\left(\int_0^T \e^{Jt}\widetilde{B} \e^{J^*t}\rd t\right)Q^*.
    \end{equation}
    
    Then, using the block decomposition in \cref{eq:Btilde_partition}, we obtain
    \begin{equation}\label{eq:W_integrand}
        \e^{Jt}\widetilde{B}\e^{J^* t}=
        \begin{pmatrix}
            \e^{J_1 t}\widetilde{B}_{1,1}\e^{J_1^* t} & \cdots & \e^{J_1 t}\widetilde{B}_{1,m}\e^{J_m^* t} \\
            \vdots & \ddots & \vdots \\
            \e^{J_m t}\widetilde{B}_{m,1}\e^{J_1^* t} & \cdots & \e^{J_m t}\widetilde{B}_{m,m}\e^{J_m^* t}
        \end{pmatrix}.
    \end{equation}
    From \cref{eq:D_definition,eq:scaledW_definition,eq:W_jordan,eq:W_integrand}, the $(i,j)$-block of $W(B;T)$ can be written as
    \begin{equation}\label{eq:scaledW_block}
        D_i(T)^{-1}\left\{\int_0^T \e^{J_i t}\widetilde{B}_{i,j}\e^{J_j^* t}\rd t\right\}D_j(T)^{-*}.
    \end{equation}
    
    In the following, we analyze the block representation~\cref{eq:scaledW_block} by distinguishing the following nine cases:
    \begin{enumerate}[label=\arabic*), ref=Case~\arabic*)]
        \item\label{case:stable_stable} $i=1,\ldots,m_-,\quad j=1,\ldots,m_-$
        \item\label{case:stable_imaginary} $i=1,\ldots,m_-,\quad j=m_-+1,\ldots,m_-+m_0$
        \item\label{case:stable_unstable} $i=1,\ldots,m_-,\quad j=m_-+m_0+1,\ldots,m$
        \item\label{case:imaginary_stable} $i=m_-+1,\ldots,m_-+m_0,\quad j=1,\ldots,m_-$
        \item\label{case:imaginary_imaginary} $i=m_-+1,\ldots,m_-+m_0,\quad j=m_-+1,\ldots,m_-+m_0$
        \item\label{case:imaginary_unstable} $i=m_-+1,\ldots,m_-+m_0,\quad j=m_-+m_0+1,\ldots,m$
        \item\label{case:unstable_stable} $i=m_-+m_0+1,\ldots,m,\quad j=1,\ldots,m_-$
        \item\label{case:unstable_imaginary} $i=m_-+m_0+1,\ldots,m,\quad j=m_-+1,\ldots,m_-+m_0$
        \item\label{case:unstable_unstable} $i=m_-+m_0+1,\ldots,m,\quad j=m_-+m_0+1,\ldots,m$
    \end{enumerate}

    \subsubsection*{\ref{case:stable_stable}}
    Noting that $D_i(T)=I_{n_i}$ and $D_j(T)=I_{n_j}$, \cref{eq:scaledW_block} reduces to
        $\int_0^T \e^{J_i t}\widetilde{B}_{i,j}\e^{J_j^* t}\rd t$.
    Since $\re(\lambda_i),\re(\lambda_j)<0$, it follows from \cref{prop:sylvester_convergence} that the $(i,j)$-block converges, as $T\to\infty$, to the unique solution $W_-^{(i,j)}(B)\in\CC^{n_i\times n_j}$ satisfying
    \begin{equation*}
        J_i W_-^{(i,j)}(B)+W_-^{(i,j)}(B)J_j^*+\widetilde{B}_{i,j}=O.
    \end{equation*}
    
    Then, the upper-left block of $W(B;T)$ of size $n_-\times n_-$ converges to
    \begin{equation*}
        \begin{pmatrix}
            W_-^{(1,1)}(B) & \cdots & W_-^{(1,m_-)}B) \\
            \vdots & \ddots & \vdots \\
            W_-^{(m_-,1)}(B) & \cdots & W_-^{(m_-,m_-)}(B)
        \end{pmatrix},
    \end{equation*}
    which coincides with $W_-(B)$.

    \subsubsection*{\ref{case:stable_imaginary}}
    Noting that $D_i(T)=I_{n_i}$, \cref{eq:scaledW_block} reduces to
    \begin{equation*}
        \left\{\int_0^T \e^{J_i t}\widetilde{B}_{i,j}\e^{J_j^* t}\rd t\right\}D_j(T)^{-*}.
    \end{equation*}
    Since $\re(\lambda_i)<0=\re(\lambda_j)$, it follows from \cref{prop:sylvester_convergence} that the integral converges as $T\to\infty$.
    Moreover, since $D_j(T)^{-*}\to0$, the $(i,j)$-block converges to zero.
    
    \subsubsection*{\ref{case:stable_unstable}}
    Noting that $D_i(T)=I_{n_i}$ and $D_j(T)=\e^{J_jT}$, \cref{eq:scaledW_block} reduces to
    \begin{equation*}
        \left\{\int_0^T \e^{J_i t}\widetilde{B}_{i,j}\e^{J_j^* t}\rd t\right\}\e^{-J_j^*T}.
    \end{equation*}
    There exists a constant $c\geq0$ such that
        $\norm{\e^{J_it}}\leq c t^{n_i-1}\e^{\re(\lambda_i)t}$,
        $\norm{\e^{J_j^*t}}\leq c t^{n_j-1}\e^{\re(\lambda_j)t}$, and 
        $\norm{\widetilde{B}_{i,j}}\leq c$
    for all $t\geq0$.
    Consequently, we obtain
    \begin{align*}
        &\norm{\int_0^T \e^{J_i t}\widetilde{B}_{i,j}\e^{J_j^* t}\rd t} \\
        &\leq\int_0^T \norm{\e^{J_i t}}\norm{\e^{J_j^* t}}\norm{\widetilde{B}_{i,j}}\rd t \\
        &\leq\int_0^T c^3 t^{n_i+n_j-2}\e^{(\re(\lambda_i)+\re(\lambda_j))t}\rd t \\
        &\leq\left\{
        \begin{aligned}
            &c' \qquad \text{if $\re(\lambda_i)+\re(\lambda_j)<0$}, \\
            &c'T^{n_i+n_j-2}\e^{(\re(\lambda_i)+\re(\lambda_j))T} \\
            &\qquad \text{if $\re(\lambda_i)+\re(\lambda_j)>0$}, \\
            &c'T^{n_i+n_j-1} \qquad \text{if $\re(\lambda_i)+\re(\lambda_j)=0$}
        \end{aligned}
        \right.
    \end{align*}
    for some constant $c'\geq0$.
    Since
    \begin{equation*}
        \norm{\e^{-J_j^*T}}\leq c T^{n_j-1}\e^{-\re(\lambda_j)T},
    \end{equation*}
    we further obtain
    \begin{align*}
        &\norm{\left\{\int_0^T \e^{J_i t}\widetilde{B}_{i,j}\e^{J_j^* t}\rd t\right\}\e^{-J_j^*T}} \\
        &\leq\norm{\int_0^T \e^{J_i t}\widetilde{B}_{i,j}\e^{J_j^* t}\rd t}\norm{\e^{-J_j^*T}} \\
        &\leq
        \begin{cases*}
            cc'T^{n_j-1}\e^{-\re(\lambda_j)T} & if $\re(\lambda_i)+\re(\lambda_j)<0$, \\
            cc'T^{n_i+2n_j-3}\e^{\re(\lambda_i)T} & if $\re(\lambda_i)+\re(\lambda_j)>0$, \\
            cc'T^{n_i+2n_j-2}\e^{-\re(\lambda_j)T} & if $\re(\lambda_i)+\re(\lambda_j)=0$.
        \end{cases*}
    \end{align*}
    Since $\re(\lambda_i)<0<\re(\lambda_j)$, the $(i,j)$-block converges to zero in all cases.
    
    \subsubsection*{\ref{case:imaginary_stable}}
    Similar to \ref{case:stable_imaginary}, the $(i,j)$-block converges to zero.
    
    \subsubsection*{\ref{case:imaginary_imaginary}}
    Since $\e^{J_it}=\sum_{k=0}^{n_i-1}\e^{\lambda_it}\dfrac{t^k}{k!}N_{n_i}^k$, we obtain
    \begin{equation*}
        \e^{J_it}\widetilde{B}_{i,j}\e^{J_j^*t}=\sum_{k=0}^{n_i-1}\sum_{\ell=0}^{n_j-1}\e^{(\lambda_i+\overline{\lambda_j})t}\dfrac{t^{k+\ell}}{k!\ell!}N_{n_i}^k\widetilde{B}_{i,j}(N_{n_j}^*)^\ell.
    \end{equation*}
    Consequently, \cref{eq:scaledW_block} can be written as
    \begin{multline}\label{eq:scaledW_block_imim}
        \sum_{k=0}^{n_i-1}\sum_{\ell=0}^{n_j-1}\left\{\int_0^T\e^{(\lambda_i+\overline{\lambda_j})t}\dfrac{t^{k+\ell}}{k!\ell!}\rd t\right\}\times \\
        D_i(T)^{-1}N_{n_i}^k\widetilde{B}_{i,j}(N_{n_j}^*)^\ell D_j(T)^{-*}.
    \end{multline}
    
    Since the potentially nonzero entries of $N_{n_i}^k\widetilde{B}_{i,j}(N_{n_j}^*)^\ell$ are confined to its upper-left $(n_i-k)\times(n_j-\ell)$ block, it follows that
    \begin{align*}
        &D_i(T)^{-1}N_{n_i}^k\widetilde{B}_{i,j}(N_{n_j}^*)^\ell D_j(T)^{-*} \\
        &=T^{-(k+\ell+1)}\widetilde{D}_{i,k}(T)N_{n_i}^k\widetilde{B}_{i,j}(N_{n_j}^*)^\ell\widetilde{D}_{j,\ell}(T)^*,
    \end{align*}
    where
    \begin{align*}
        \widetilde{D}_{i,k}(T)&\coloneq\diag(T^{-(n_i-k-1)},\ldots,T^{-1},1,0,\ldots,0), \\
        \widetilde{D}_{j,\ell}(T)&\coloneq\diag(T^{-(n_j-\ell-1)},\ldots,T^{-1},1,0,\ldots,0).
    \end{align*}
    We define $E_{i,k}\coloneq\lim_{T\to\infty} \widetilde{D}_{i,k}(T)$, where $E_{i,k}$ has a $1$ in the $(n_i-k,n_i-k)$-entry and zeros elsewhere; similarly, we define $E_{j,\ell}\coloneq\lim_{T\to\infty} \widetilde{D}_{j,\ell}(T)$.
    
    Moreover, noting that $\re(\lambda_i)=\re(\lambda_j)=0$, we obtain
    \begin{equation*}
        \left\{
        \begin{aligned}
            \left|\int_0^T\e^{(\lambda_i+\overline{\lambda_j})t}\dfrac{t^{k+\ell}}{k!\ell!}\rd t\right|&\leq cT^{k+\ell} & \text{if $\lambda_i\neq\lambda_j$}, \\
            \int_0^T\e^{(\lambda_i+\overline{\lambda_j})t}\dfrac{t^{k+\ell}}{k!\ell!}\rd t&=\dfrac{T^{k+\ell+1}}{k!\ell!(k+\ell+1)} & \text{if $\lambda_i=\lambda_j$}
        \end{aligned}
        \right.
    \end{equation*}
    for some constant $c\geq 0$.
    Consequently, the summand of \cref{eq:scaledW_block_imim} converges, as $T\to\infty$, to
    \begin{equation*}
        \left\{
        \begin{aligned}
            &O & \text{if $\lambda_i\neq\lambda_j$}, \\
            &\dfrac{1}{k!\ell!(k+\ell+1)}E_{i,k}N_{n_i}^k\widetilde{B}_{i,j}(N_{n_j}^*)^\ell E_{j,\ell}^* & \text{if $\lambda_i=\lambda_j$}.
        \end{aligned}
        \right.
    \end{equation*}
    Then, $E_{i,k}N_{n_i}^k\widetilde{B}_{i,j}(N_{n_k}^*)^\ell E_{j,\ell}^*\in\CC^{n_i\times n_j}$ has a $\widetilde{b}_{i,j}$ in the $(n_i-k,n_j-\ell)$-entry and zeros elsewhere.
    Hence, as $T\to\infty$, the $(i,j)$-block converges to $\widetilde{b}_{i,j}C_{i,j}$, where $C_{i,j}$ is defined in \cref{eq:C_definition}, when $\lambda_i=\lambda_j$, whereas it converges to zero when $\lambda_i\neq\lambda_j$.
    Furthermore, the central block of $W(B;T)$ of size $n_0\times n_0$ converges to $W_0(B)$.

    \subsubsection*{\ref{case:imaginary_unstable}}
    Noting that $D_j(T)=\e^{J_j T}$, we can rewrite \cref{eq:scaledW_block} as
    \begin{align*}
        &D_i(T)^{-1}\e^{J_iT}\left\{\int_0^T \e^{J_i(t-T)}\widetilde{B}_{i,j}\e^{J_j^*(t-T)}\rd t\right\} \\
        &=D_i(T)^{-1}\e^{J_iT}\left\{\int_0^T \e^{-J_i t}\widetilde{B}_{i,j}\e^{-J_j^*t}\rd t\right\}.
    \end{align*}
    Since $\re(\lambda_i)=0<\re(\lambda_j)$, it follows from \cref{prop:sylvester_convergence} that the integral of the right-hand side converges as $T\to\infty$.
    Moreover, a direct calculation yields
    \begin{align*}
    &D_i(T)^{-1}\e^{J_iT} \\
    &=
    \e^{\lambda_i T}
    \begin{pmatrix}
        T^{-(n_i-\frac{1}{2})}
        & T^{-(n_i-\frac{3}{2})}
        & \cdots
        & \dfrac{1}{(n_i-1)!}T^{-\frac{1}{2}}
        \\
        0
        & T^{-(n_i-\frac{3}{2})}
        & \cdots
        & \dfrac{1}{(n_i-2)!}T^{-\frac{1}{2}}
        \\
        \vdots
        & \ddots
        & \ddots
        & \vdots
        \\
        0
        & \cdots
        & 0
        & T^{-\frac{1}{2}}
    \end{pmatrix},
\end{align*}
    and thus $D_i(T)^{-1}\e^{J_iT}\to0$ as $T\to\infty$.
    Hence, the $(i,j)$-block converges to zero.
    
    \subsubsection*{\ref{case:unstable_stable}}
    Similar to \ref{case:stable_unstable}, the $(i,j)$-block converges to zero.
    
    \subsubsection*{\ref{case:unstable_imaginary}}
    Similar to \ref{case:imaginary_unstable}, the $(i,j)$-block converges to zero.
    
    \subsubsection*{\ref{case:unstable_unstable}}
    Noting that $D_i(T)=\e^{J_iT}$ and $D_j(T)=\e^{J_jT}$, \cref{eq:scaledW_block} reduces to
        $\int_0^T \e^{J_i(t-T)}\widetilde{B}_{i,j}\e^{J_j^*(t-T)}\rd t
        =\int_0^T \e^{-J_it}\widetilde{B}_{i,j}\e^{-J_j^*t}\rd t$.
    Since $\re(\lambda_i),\re(\lambda_j)>0$, it follows from \cref{prop:sylvester_convergence} that the integral of the right-hand side converges, as $T\to\infty$, to the unique solution $W_+^{(i,j)}(B)\in\CC^{n_i\times n_j}$ satisfying
    \begin{equation*}
        (-J_i)W_+^{(i,j)}(B)+W_+^{(i,j)}(B)(-J_j^*)+\widetilde{B}_{i,j}=0.
    \end{equation*}
    
    Then, the lower-right block of $W(B;T)$ of size $n_+\times n_+$ converges to
    \begin{equation*}
        \begin{pmatrix}
            W_+^{(\gamma+1,\gamma+1)}(B) & \cdots & W_+^{(\gamma+1,m)}(B) \\
            \vdots & \ddots & \vdots \\
            W_+^{(m,\gamma+1)}(B) & \cdots & W_+^{(m,m)}(B)
        \end{pmatrix},
    \end{equation*}
    where we set $\gamma\coloneq m_-+m_0$ to simplify the indexing.
    The matrix coincides with $W_+(B)$.

    Summarizing the discussion so far, $W(B;T)$ converges to $W_\infty(B)$, completing the proof.
    \qed
\end{proof}

\section{Proof of \cref{prop:finite_nonempty_lemma,prop:infinite_nonempty_lemma}}
\label{sec:proof_nonempty_lemma}

First, we prove \cref{prop:finite_nonempty_lemma}.

\begin{proof}[\cref{prop:finite_nonempty_lemma}]
    Let $y\in\RR^n$ be a nonzero vector.
    Since $\e^{A^\top t}$ is nonsingular for all $t$, $\e^{A^\top t}y\neq0$.
    Therefore, since $\diag(p)$ is positive definite,
    \begin{align*}
        y^\top\widetilde{W}(p;T)y=\int_0^T y^\top\e^{At}\diag(p)\e^{A^\top t}y \rd t >0.
    \end{align*}
    This proves that $\widetilde{W}(p;T)$ is positive definite.
    \qed
\end{proof}

Next, to prove \cref{prop:infinite_nonempty_lemma}, we employ the following lemma.

\begin{lemma}[{\cite[Theorem 12]{Lancaster1970}}]\label{prop:lyapunov_lemma}
    Let $C\in\CC^{\ell\times\ell}$ be Hurwitz and $P\in\CC^{\ell\times\ell}$ be Hermitian positive definite.
    Then, the Lyapunov equation
        $CY+YC^*+P=O$
    admits a unique solution $Y\in\CC^{\ell\times\ell}$, which is Hermitian positive definite.
\end{lemma}

\begin{proof}[\cref{prop:infinite_nonempty_lemma}]
    From \cref{eq:W_block_diagonal}, it suffices to prove that $W_-(p), W_0(p)$, and $W_+(p)$ are positive definite.
    Let $B\coloneq\diag(p), \widetilde{B}\coloneq Q^{-1}BQ^{-*}$, and partition $\widetilde{B}$ in the same way as \cref{eq:Btilde_partition}.
    Then, the matrices $\widetilde{B}_{-,-}, \widetilde{B}_{0,0}$, and $\widetilde{B}_{+,+}$ are positive definite since they are principal submatrices of $\widetilde{B}$, which is positive definite.

    From \cref{eq:W_stable_jordan_sylvester}, the matrix $W_-(p)$ is the unique solution to
    \begin{equation*}
        J_-W_-(p)+W_-(p)J_-^*+\widetilde{B}_{-,-}=O,
    \end{equation*}
    and it follows from \cref{prop:lyapunov_lemma} that $W_-(p)$ is positive definite.
    Similarly, we can prove that $W_+(p)$ is positive definite.

    We next prove that $W_0(p)$ is positive definite.
    Since $W_0(p)$ is a block-diagonal matrix from \cref{eq:W0_block_diagonal}, it suffices to show that $W_0^{(\ell)}(p)$ is positive definite for $\ell=1,\ldots,m_0'$.
    To simplify the indexing, let $\alpha=m_-+m_0^{(1)}+\cdots+m_0^{(\ell-1)},\ \beta=\alpha+m_0^{(\ell)}$.
    Let $e^{(i)}\coloneq(0,\ldots,0,1)^\top\in\RR^{n_i}$.
    Then, a direct calculation yields
    \begin{equation*}
        C_{i,j}=\int_0^1 \e^{N_i t}e^{(i)}e^{(j)}{}^\top\e^{N_j^\top t}\rd t.
    \end{equation*}

    Furthermore, let
    \begin{equation*}
        \widetilde{B}_{0,0}^{(\ell)}\coloneq
        \begin{pmatrix}
            \widetilde{b}_{\alpha+1,\alpha+1} & \cdots & \widetilde{b}_{\alpha+1,\beta} \\
            \vdots & \ddots & \vdots \\
            \widetilde{b}_{\beta,\alpha+1} & \cdots & \widetilde{b}_{\beta,\beta}
        \end{pmatrix}
        \in\CC^{m_0^{(\ell)}\times m_0^{(\ell)}},
    \end{equation*}
    which is positive definite since it is a principal submatrix of $\widetilde{B}_{0,0}$,
    \begin{equation*}
        E^{(\ell)}\coloneq
        \begin{pmatrix}
            e^{(\alpha+1)} & \cdots & 0 \\
            \vdots & \ddots & \vdots \\
            0 & \cdots & e^{(\beta)}
        \end{pmatrix}\in\RR^{n_0^{(\ell)}\times m_0^{(\ell)}},
    \end{equation*}
    and
    \begin{equation*}
        J_0^{(\ell)}\coloneq \mathrm{diag}(J_{\alpha+1},\ldots, J_\beta)
        \in\CC^{n_0^{(\ell)}\times n_0^{(\ell)}}.
    \end{equation*}
    By further calculation, it follows from \cref{eq:W0l_block} that
    \begin{equation*}
        W_0^{(\ell)}(p)=\int_0^1 \e^{J_0^{(\ell)}t}E^{(\ell)}\widetilde{B}_{0,0}^{(\ell)}E^{(\ell)}{}^\top \e^{J_0^{(\ell)}{}^* t}\rd t.
    \end{equation*}

    Let $y\in\CC^{n_0^{(\ell)}}$ be a nonzero vector, and partition it into blocks as
        $y
        =
        \left(
            y^{(\alpha+1)}{}^\top,
            \ldots,
            y^{(\beta)}{}^\top
        \right)^\top$.
    Then,
    \begin{align*}
        y^*W_0^{(\ell)}(p)y
        =
        \int_0^1
        y^*\e^{J_0^{(\ell)}t}
        E^{(\ell)}
        \widetilde{B}_{0,0}^{(\ell)}
        E^{(\ell)}{}^\top
        \e^{J_0^{(\ell)}{}^*t}y
        \rd t.
    \end{align*}
    Moreover,
    \begin{align*}
        y^*\e^{J_0^{(\ell)}t}E^{(\ell)}
        =
        \bigl(
            y^{(\alpha+1)}{}^*
            \e^{J_{\alpha+1}t}e^{(\alpha+1)},
            \ldots,
            y^{(\beta)}{}^*
            \e^{J_\beta t}e^{(\beta)}
        \bigr).
    \end{align*}

    Since $y$ is nonzero, there exists
    $i\in\{\alpha+1,\ldots,\beta\}$ such that $y^{(i)}\neq0$.
    Recalling that $J_i=\lambda_i I_{n_i}+N_i$, we have
    \begin{align*}
        y^{(i)}{}^*\e^{J_it}e^{(i)}
        =
        \e^{\lambda_i t}
        y^{(i)}{}^*\e^{N_it}e^{(i)}
        =
        \e^{\lambda_i t}
        \sum_{k=0}^{n_i-1}
        \frac{\overline{y^{(i)}_{n_i-k}}}{k!}t^k.
    \end{align*}
    The polynomial on the right-hand side is not identically zero
    because $y^{(i)}\neq0$.
    Furthermore, $\e^{\lambda_i t}$ never vanishes.
    Consequently,
    $y^{(i)}{}^*\e^{J_it}e^{(i)}$ can vanish at only finitely many
    points in $[0,1]$.
    It follows that
        $y^*\e^{J_0^{(\ell)}t}E^{(\ell)}
        \neq 0$
    for almost every $t\in[0,1]$.

    Since $\widetilde{B}_{0,0}^{(\ell)}$ is positive definite, we have
        $y^*\e^{J_0^{(\ell)}t}
        E^{(\ell)}
        \widetilde{B}_{0,0}^{(\ell)}
        E^{(\ell)}{}^\top
        \e^{J_0^{(\ell)}{}^*t}y
        >0$
    for almost every $t\in[0,1]$.
    Therefore,
        $y^*W_0^{(\ell)}(p)y>0$.
    Since $y$ was arbitrary, $W_0^{(\ell)}(p)$ is positive definite.
    This completes the proof. \qed
\end{proof}

\section{Proof of \cref{prop:vcs_convergence,prop:aecs_convergence}}
\label{sec:convergence_proof}

Before proving \cref{prop:vcs_convergence,prop:aecs_convergence}, we collect several definitions and preliminary facts that will be used throughout the proofs.

Let $\widetilde{\Delta}$ be the simplex in $\RR^{n-1}$, i.e.,
\begin{equation*}
    \widetilde{\Delta}\coloneq
    \left\{
        (\widetilde{p}_i)\in\RR^{n-1}
    \,\middle|\,
        \begin{gathered}
            \widetilde{p}_i\geq0 \ (i=1,\ldots,n-1), \\
            \textstyle \sum_{i=1}^{n-1} \widetilde{p}_i \leq 1
        \end{gathered}
    \right\}.
\end{equation*}
The mapping $\varphi\colon\widetilde{\Delta}\to\Delta$ defined by
\begin{equation*}
    \textstyle \varphi(\widetilde{p}_1,\ldots,\widetilde{p}_{n-1})=(\widetilde{p}_1,\ldots,\widetilde{p}_{n-1},1-\sum_{i=1}^{n-1}\widetilde{p}_i)^\top
\end{equation*}
is a homeomorphism; its inverse is given by
\begin{equation*}
    \varphi^{-1}(p_1,\ldots,p_n)=(p_1,\ldots,p_{n-1})^\top.
\end{equation*}
For a set $\Omega\subset\RR^{n-1}$, $\interior\Omega$ denotes its interior.

Let $f\colon\RR^{n-1}\to\RR\cup\{\infty\}$.
Define $\dom f\coloneq\{\widetilde{p}\in\RR^{n-1}\mid f(\widetilde{p})<\infty\}$.
We say that $f$ is proper if $\dom f\neq\emptyset$.
We say that $f$ is lower semicontinuous if the sublevel set $\mathrm{lev}_{\leq\alpha}f\coloneq\{\widetilde{p}\in\RR^{n-1}\mid f(\widetilde{p})\leq\alpha\}$ is closed for all $\alpha\in\RR$.
We say that $f$ is level-bounded if $\mathrm{lev}_{\leq\alpha}f$ is bounded for all $\alpha\in\RR$.

Let $\{f^{(k)}\}_{k\in\NN}$ be a sequence of functions $f^{(k)}\colon\RR^{n-1}\to\RR\cup\{\infty\}$.
We say that $\{f^{(k)}\}_{k\in\NN}$ is equi level-bounded if for every $\alpha\in \RR$, there exists a bounded set $B_{\alpha}\subset\RR^{n-1}$ such that $\mathrm{lev}_{\leq\alpha}f^{(k)}\subset B_{\alpha}$ for all $k$.
We say that $\{f^{(k)}\}_{k\in\NN}$ epi-converges to $f$ if, for every $\widetilde{p}\in\RR^{n-1}$,
\begin{equation*}
    \left\{
        \begin{gathered}
        \liminf_{k\to\infty} f^{(k)}(\widetilde{p}^{(k)})\geq f(\widetilde{p}) \quad \text{for every sequence $\widetilde{p}^{(k)}\to \widetilde{p}$,} \\
        \limsup_{k\to\infty} f^{(k)}(\widetilde{p}^{(k)})\leq f(\widetilde{p}) \quad \text{for some sequence $\widetilde{p}^{(k)}\to \widetilde{p}$}.
        \end{gathered}
    \right.
\end{equation*}

We employ the following three propositions.
\begin{proposition}[{\cite[Theorem 1.9]{RockafellarWets1998}}]\label{prop:attainment_lemma}
    Assume that $f\colon\RR^{n-1}\to\RR\cup\{\infty\}$ is proper, lower semicontinuous, and level-bounded.
    Then, $\argmin f$ is nonempty and compact, and the minimum value of $f$ is finite.
\end{proposition}

\begin{proposition}[{\cite[Theorem 7.17]{RockafellarWets1998}}]\label{prop:convergence_lemma}
    Assume that $\{f^{(k)}\}_{k\in\NN}$ is a sequence of convex functions $f^{(k)}\colon\RR^{n-1}\to\RR\cup\{\infty\}$, and that $f\colon\RR^{n-1}\to\RR\cup\{\infty\}$ is a convex and lower semicontinuous function such that $\interior\dom f\neq\emptyset$.
    If there is a dense subset $\Omega$ of $\RR^{n-1}$ such that $f^{(k)}(\widetilde{p})\to f(\widetilde{p})$ for all $\widetilde{p}\in \Omega$, then $\{f^{(k)}\}_{k\in\NN}$ epi-converges to $f$.
\end{proposition}

\begin{proposition}[{\cite[Theorem 7.33]{RockafellarWets1998}}]\label{prop:epi_optimal}
    Assume that
    $f^{(k)}\colon\RR^{n-1}\to\RR\cup\{\infty\}$
    and $f\colon\RR^{n-1}\to\RR\cup\{\infty\}$
    are proper lower semicontinuous functions.
    Suppose that
    $\{f^{(k)}\}_{k\in\NN}$ is equi level-bounded and epi-converges to $f$.
    Then, any sequence $\{\widetilde{p}^{(k)}\}_{k\in\NN}$ satisfying $\widetilde{p}^{(k)}\in\argmin f^{(k)}$ is bounded, and all its cluster points belong to $\argmin f$.
    In particular, if $\argmin f=\{\widetilde{p}^*\}$, then $\widetilde{p}^{(k)} \to \widetilde{p}^*$.
\end{proposition}

Furthermore, we prepare the following lemma to prove \cref{prop:aecs_convergence}.
\begin{lemma}\label{prop:trace_inequality}
    Let $C_1, C_2\in\CC^{\ell\times\ell}$ be Hermitian positive definite matrices, and let $\mu>0$ be the minimum eigenvalue of $C_2$.
    Then,
        $\tr(C_1 C_2)\geq\mu\tr(C_1)$
    holds.
\end{lemma}

\begin{proof}
    Since $C_2\succeq\mu I_\ell$, we obtain $C_1^{1/2} C_2 C_1^{1/2}\succeq\mu C_1$.
    Therefore,
        $\tr(C_1 C_2)=\tr(C_1^{1/2}C_2C_1^{1/2})\geq\mu\tr(C_1)$.
    \qed
\end{proof}

Then, we can prove \cref{prop:vcs_convergence,prop:aecs_convergence}.

\begin{proof}[\cref{prop:vcs_convergence}]
    Let $\{T_k\}_{k\in\NN}$ be any sequence such that $T_k>0$ and $T_k\to\infty$ as $k\to\infty$.
    It suffices to prove that $p_{\mathrm{VCS}}^{T_k}\to p_{\mathrm{VCS}}^{\infty}$ as $k\to\infty$.

    Define $g^{(k)}\colon\RR^{n-1}\to\RR\cup\{\infty\}$ by
    \begin{equation*}
        g^{(k)}(\widetilde{p})\coloneq\left\{
        \begin{aligned}
            &-\log\det W(\varphi(\widetilde{p});T_k) \\
            &\qquad \text{if $\widetilde{p}\in\widetilde{\Delta},\ W(\varphi(\widetilde{p});T_k)\succ O$,} \\
            & \infty \quad \text{otherwise},
        \end{aligned}
        \right.
    \end{equation*}
    and define $g\colon\RR^{n-1}\to\RR\cup\{\infty\}$ by
    \begin{equation*}
        g(\widetilde{p})\coloneq\left\{
        \begin{aligned}
            &-\log\det W_\infty(\varphi(\widetilde{p})) \\
            &\qquad \text{if $\widetilde{p}\in\widetilde{\Delta},\ W_\infty(\varphi(\widetilde{p}))\succ O$,} \\
            & \infty \quad \text{otherwise}.
        \end{aligned}
        \right.
    \end{equation*}
    Then, $\{g^{(k)}\}_{k\in\NN}$ is a sequence of convex functions and equi level-bounded since $\mathrm{lev}_{\leq\alpha}g^{(k)}\subset\widetilde{\Delta}$ for every $\alpha\in\RR$.
    Similarly, $g$ is convex and level-bounded.
    It follows from \cref{prop:finite_nonempty_lemma} that $g^{(k)}$ is proper.
    Furthermore, it follows from \cref{prop:infinite_nonempty_lemma} that $\interior\dom g=\interior\widetilde{\Delta}\neq\emptyset$.

    Next, we prove that $g$ is lower semicontinuous, that is, $\mathrm{lev}_{\leq\alpha}g$ is closed for every $\alpha\in\RR$.
    Let $\widetilde{p}^{(k)}\in\mathrm{lev}_{\leq\alpha}g$ with $\widetilde{p}^{(k)}\to\widetilde{p}\in\RR^{n-1}$ as $k\to\infty$.
    Since $\det W_{\infty}(\varphi(\widetilde{p}^{(k)}))\geq \e^{-\alpha}$, by continuity, we obtain $\det W_{\infty}(\varphi(\widetilde{p}))\geq \e^{-\alpha}$ and $W_{\infty}(\varphi(\widetilde{p}))\succ O$.
    Since $\widetilde{\Delta}$ is closed, we obtain $\widetilde{p}\in \widetilde{\Delta}$; hence $g(\widetilde{p})\leq\alpha$, implying $\widetilde{p}\in\mathrm{lev}_{\leq\alpha}g$.
    This shows that $\mathrm{lev}_{\leq\alpha}g$ is closed.
    Similarly, we can prove that $g^{(k)}$ is lower semicontinuous.

    Thus, it follows from \cref{prop:attainment_lemma} that $\emptyset\neq\argmin g^{(k)}\subset\widetilde{\Delta}$ and $\emptyset\neq\argmin g\subset\widetilde{\Delta}$.
    If $\widetilde{p}^{(k)}\in\argmin g^{(k)}$, then $\varphi(\widetilde{p}^{(k)})$ is an optimal solution to \cref{prob:scaled_vcs}, and hence also to \cref{prob:vcs}.
    If $\widetilde{p}^*\in\argmin g$, then $\varphi(\widetilde{p}^*)$ is an optimal solution to \cref{prob:infinite_vcs}.
    Furthermore, by assumption, $\argmin g=\{\widetilde{p}^*\}$.
    Therefore, by continuity of $\varphi$, it suffices to prove that $\widetilde{p}^{(k)}\to\widetilde{p}^*$ as $k\to\infty$.

    Let $\Omega\coloneq\interior\dom\widetilde{\Delta}\cup(\RR^{n-1}\setminus\widetilde{\Delta})$, which is a dense subset of $\RR^{n-1}$.
    If $\widetilde{p}\in\RR^{n-1}\setminus\widetilde{\Delta}$, then $g^{(k)}(\widetilde{p})=\infty\to g(\widetilde{p})=\infty$.
    If $\widetilde{p}\in\interior\dom\widetilde{\Delta}$, then it follows from \cref{prop:finite_nonempty_lemma,prop:infinite_nonempty_lemma} that $W(\varphi(\widetilde{p});T_k)\succ O$ and $W_\infty(\varphi(\widetilde{p}))\succ O$.
    By continuity, we obtain $g^{(k)}(\widetilde{p})=-\log\det W(\varphi(\widetilde{p});T_k)\to g(\widetilde{p})=-\log\det W_\infty(\varphi(\widetilde{p}))$.
    Consequently, $g^{(k)}(\widetilde{p})\to g(\widetilde{p})$ for all $\widetilde{p}\in \Omega$, and it follows from \cref{prop:convergence_lemma} that $\{g^{(k)}\}_{k\in\NN}$ epi-converges to $g$.
    It follows from \cref{prop:epi_optimal} that $\widetilde{p}^{(k)}\to\widetilde{p}^*$, and this completes the proof.
    \qed
\end{proof}

\begin{proof}[\cref{prop:aecs_convergence}]
    Let $\{T_k\}_{k\in\NN}$ be any sequence such that $T_k>0$ and $T_k\to\infty$ as $k\to\infty$.
    It suffices to prove that $p_{\mathrm{AECS}}^{T_k}\to p_{\mathrm{AECS}}^{\infty}$ as $k\to\infty$.

    Define $h^{(k)}\colon\RR^{n-1}\to\RR\cup\{\infty\}$ by
    \begin{equation*}
        h^{(k)}(\widetilde{p})\coloneq\left\{
        \begin{aligned}
            &\tr\left(Q^{-*}D(T_k)^{-*}W(\varphi(\widetilde{p});T_k)^{-1}D(T_k)^{-1}Q^{-1}\right) \\
            &\qquad \text{if $\widetilde{p}\in\widetilde{\Delta},\ W(\varphi(\widetilde{p});T_k)\succ O$,} \\
            & \infty \quad \text{otherwise},
        \end{aligned}
        \right.
    \end{equation*}
    and define $h\colon\RR^{n-1}\to\RR\cup\{\infty\}$ by
    \begin{equation*}
        h(\widetilde{p})\coloneq\left\{
        \begin{aligned}
            &\tr\left(Q^{-*} \mathrm{diag}(W_-(\varphi(\widetilde{p}))^{-1}, O, O) Q^{-1}\right) \\
            &\qquad \text{if $\widetilde{p}\in\widetilde{\Delta},\ W_-(\varphi(\widetilde{p}))\succ O$,} \\
            & \infty \quad \text{otherwise}.
        \end{aligned}
        \right.
    \end{equation*}
    Then, $\{h^{(k)}\}_{k\in\NN}$ is a sequence of convex functions and equi level-bounded since $\mathrm{lev}_{\leq\alpha}h^{(k)}\subset\widetilde{\Delta}$ for every $\alpha\in\RR$.
    Similarly, $h$ is convex and level-bounded.
    It follows from \cref{prop:finite_nonempty_lemma} that $h^{(k)}$ is proper.
    Furthermore, it follows from \cref{prop:infinite_nonempty_lemma} that $\interior\dom h=\interior\widetilde{\Delta}\neq\emptyset$.

    Next, we prove that $h$ is lower semicontinuous.
    Let $(Q^{-1}Q^{-*})_{-}$ be the upper-left block of $Q^{-1}Q^{-*}$ of size $n_-\times n_-$.
    Then, we obtain
    \begin{align*}
        &\tr\left(Q^{-*}
            \mathrm{diag}(W_-(\varphi(\widetilde{p}))^{-1}, O, O)
            Q^{-1}\right) \\
        &=\tr\left(W_-(\varphi(\widetilde{p}))^{-1}(Q^{-1}Q^{-*})_{-}\right).
    \end{align*}
    Furthermore, since $Q^{-1}Q^{-*}$ is Hermitian positive definite, $(Q^{-1}Q^{-*})_{-}$ is also Hermitian positive definite, and let $\mu>0$ denote its minimum eigenvalue.
    
    Let $\alpha\in\RR$ be arbitrary, and let $\widetilde{p}^{(k)}\in\mathrm{lev}_{\leq\alpha}h$ with $\widetilde{p}^{(k)}\to\widetilde{p}\in\RR^{n-1}$ as $k\to\infty$.
    Since $\widetilde{\Delta}$ is closed, we obtain $\widetilde{p}\in\Delta$.
    Suppose that $W_-(\varphi(\widetilde{p}))\not\succ O$.
    Then, $\tr\left(W_-(\varphi(\widetilde{p}^{(k)}))^{-1}\right)\to\infty$.
    It follows from \cref{prop:trace_inequality} that
    \begin{align*}
        h(\widetilde{p}^{(k)})&=\tr\left(W_-(\varphi(\widetilde{p}^{(k)}))^{-1}(Q^{-1}Q^{-*})_{-}\right) \\
        &\geq \mu\tr\left(W_-(\varphi(\widetilde{p}^{(k)}))^{-1}\right)\to\infty.
    \end{align*}
    However, this contradicts that $h(\widetilde{p}^{(k)})\leq\alpha$ for all $k$.
    Hence, $W_-(\varphi(\widetilde{p}))\succ O$.
    
    By continuity, we obtain
    \begin{equation*}
        h(\widetilde{p})=\tr\left(W_-(\varphi(\widetilde{p}))^{-1}(Q^{-1}Q^{-*})_{-}\right)\leq\alpha,
    \end{equation*}
    and therefore $\widetilde{p}\in\mathrm{lev}_{\leq\alpha}h$.
    This shows that $\mathrm{lev}_{\leq\alpha}h$ is closed, implying that $h$ is lower semicontinuous.
    Similarly, we can prove that $h^{(k)}$ is lower semicontinuous.

To apply Proposition~\ref{prop:convergence_lemma}, we verify pointwise convergence on the dense
set
\begin{equation*}
    \Omega
    \coloneq
    \operatorname{int}\widetilde{\Delta}
    \cup
    \left(\mathbb{R}^{n-1}\setminus\widetilde{\Delta}\right).
\end{equation*}
For
$\widetilde p\in\mathbb{R}^{n-1}\setminus\widetilde{\Delta}$,
the convergence is immediate because
$h^{(k)}(\widetilde p)=h(\widetilde p)=\infty$.
It therefore remains to establish the convergence for
$\widetilde p\in\operatorname{int}\widetilde{\Delta}$.

Let $\widetilde p\in\operatorname{int}\widetilde{\Delta}$ and set
$p=\phi(\widetilde p)$.
Then $p_i>0$ for every $i$, and hence Lemma~\ref{prop:infinite_nonempty_lemma} gives
    $W_\infty(p)\succ O$.
Therefore, by \eqref{eq:aecs_pointwise_matrix_limit},
    $\lim_{k\to \infty} h^{(k)}(\widetilde p)
    = h(\widetilde p)$ for every $\widetilde p\in\Omega$.
Proposition~\ref{prop:convergence_lemma} now implies that
$\{h^{(k)}\}_{k\in\mathbb N}$ epi-converges to $h$.
The conclusion then follows from \cref{prop:epi_optimal} in the same manner as
in the proof of \cref{prop:vcs_convergence}.
    \qed
\end{proof}

\begin{proof}[\cref{cor:nonunique_convergence}]
    We give the VCS argument; the AECS case is identical under Assumption~\ref{assump:stable_eigenvalue}.  For any sequence $T_k\to\infty$, use the functions $g^{(k)}$ and $g$ constructed in the proof of Theorem~\ref{prop:vcs_convergence}.  They are proper and lower semicontinuous, $\{g^{(k)}\}$ is equi level-bounded, and $g^{(k)}$ epi-converges to $g$.  Proposition~\ref{prop:epi_optimal} therefore shows that every cluster point of any sequence $\widetilde p^{(k)}\in\argmin g^{(k)}$ belongs to $\argmin g$.  The homeomorphism $\varphi$ gives the first assertion for the original allocations.

    To prove convergence to the solution set, suppose to the contrary that there are $\varepsilon>0$ and $T_k\to\infty$ such that
        $\inf_{q\in\mathcal P_{\mathrm{VCS}}^\infty}
        \norm{p_{\mathrm{VCS}}^{T_k}-q}\geq\varepsilon$.
    The simplex $\Delta$ is compact, so a subsequence of $p_{\mathrm{VCS}}^{T_k}$ converges.  Its limit belongs to $\mathcal P_{\mathrm{VCS}}^\infty$ by the first assertion, contradicting the inequality.  Hence the distance converges to zero.
    \qed
\end{proof}

\section{Proof of \cref{prop:infinite_equivalence}}
\label{sec:infinite_equivalence_proof}

\begin{proof}[\cref{prop:infinite_equivalence}]
    Because $A$ is block-diagonalized as in \cref{eq:A_block_diagonal}, there are invertible matrices $S_-$, $S_0$, and $S_+$ satisfying
    \begin{equation*}
        S_-^{-1}A_-S_-=J_-,\qquad
        S_0^{-1}A_0S_0=J_0,\qquad
        S_+^{-1}A_+S_+=J_+.
    \end{equation*}
    With $S\coloneq\diag(S_-,S_0,S_+)$, the matrices can be chosen so that $RS=Q$.  Hence
    \begin{equation*}
        \begin{pmatrix}q_{-,i}\\q_{0,i}\\q_{+,i}\end{pmatrix}
        =Q^{-1}e_i
        =S^{-1}\begin{pmatrix}r_{-,i}\\r_{0,i}\\r_{+,i}\end{pmatrix},
    \end{equation*}
    and thus $r_{\alpha,i}=S_\alpha q_{\alpha,i}$ for $\alpha\in\{-,0,+\}$.

    Let $W_{-,i}\coloneq W_-(e_ie_i^\top)$ and $W_{+,i}\coloneq W_+(e_ie_i^\top)$.  From \cref{eq:W_jordan_sylvester}, these are the unique solutions of
    \begin{subequations}\label{eq:W_jordan_sylvester_individual}
    \begin{align}
        J_-W_{-,i}+W_{-,i}J_-^*+q_{-,i}q_{-,i}^*&=O,\label{eq:W_stable_jordan_sylvester2}\\
        (-J_+)W_{+,i}+W_{+,i}(-J_+^*)+q_{+,i}q_{+,i}^*&=O.\label{eq:W_unstable_jordan_sylvester2}
    \end{align}
    \end{subequations}
    Multiplying these equations by $S_\alpha$ and $S_\alpha^*$ and using uniqueness of the corresponding Lyapunov solutions gives
    \begin{equation*}
        V_{-,i}=S_-W_{-,i}S_-^*,\qquad
        V_{+,i}=S_+W_{+,i}S_+^*.
    \end{equation*}
    Under Assumption~\ref{assump:imaginary_eigenvalue}, the center block is semisimple at zero, so $W_{0,i}=q_{0,i}q_{0,i}^*$ and $V_{0,i}=S_0W_{0,i}S_0^*$.  By linearity,
    \begin{equation}\label{eq:W_V_congruence}
        W_\infty(p)=S^{-1}V(p)S^{-*}.
    \end{equation}

    Congruence preserves positive definiteness, and
    \begin{equation*}
        -\log\det W_\infty(p)
        =-\log\det V(p)-\log\det(S^{-1}S^{-*}).
    \end{equation*}
    The last term is independent of $p$.  Therefore, the VCS feasible sets and minimizers in \cref{prob:infinite_vcs,prob:schur_vcs} coincide.

    It remains to compare the AECS objectives.  Equation~\eqref{eq:W_V_congruence} yields
    $\diag(W_-(p)^{-1},O,O)
    =
    S^*
    \diag (V_-(p)^{-1},O,O )
    S$.
    Since $SQ^{-1}=R^{-1}=(R^{(2)})^{-1}(R^{(1)})^\top$, $R^{(1)}$ is orthogonal, and $(R^{(2)})^{-1}$ is block lower triangular with an identity leading diagonal block, cyclic invariance of the trace gives
    \begin{align*}
        &\tr\!\left(Q^{-*}
        \diag(W_-(p)^{-1},O,O)Q^{-1}\right)\\
        &=\tr\!\left((R^{(2)})^{-\top}
        \diag(V_-(p)^{-1},O,O)(R^{(2)})^{-1}\right)\\
        &=\tr\!\left(V_-(p)^{-1}\right).
    \end{align*}
    The stable-block congruence also shows that $W_-(p)\succ O$ if and only if $V_-(p)\succ O$.  Thus, under Assumption~\ref{assump:stable_eigenvalue}, the AECS feasible sets and objective values coincide as claimed.
    \qed
\end{proof}

\section{Proof of \cref{prop:infinite_vcs_uniqueness,prop:infinite_aecs_uniqueness}}
\label{sec:infinite_uniqueness_proof}

To prove \cref{prop:infinite_vcs_uniqueness,prop:infinite_aecs_uniqueness}, we prepare the following lemmas, which can be proved in the same manner as \cite[Theorems 1 and 3]{SatoTerasaki2024}, although \cref{prob:aecs,prob:infinite_aecs} have slightly different forms.

\begin{lemma}\label{prop:vcs_linearly_independent}
    Assume that $W_{\infty,1},\ldots,W_{\infty,n}$ are linearly independent over $\RR$.
    Then, \cref{prob:infinite_vcs} admits a unique optimal solution.
\end{lemma}

\begin{lemma}\label{prop:aecs_linearly_independent}
    Assume that $W_{-,1},\ldots,W_{-,n}$ are linearly independent over $\RR$.
    Then, \cref{prob:infinite_aecs} admits a unique optimal solution.
\end{lemma}

\begin{proof}[\cref{prop:infinite_vcs_uniqueness}]
    Let $x=(x_1,\ldots,x_n)^\top\in\RR^n$ be any vector such that $W_\infty(x)=0$.
    Then, from \cref{prop:vcs_linearly_independent}, it suffices to prove $x=0$.

    Let $X\coloneq\diag(x)$ and $\widetilde{X}\coloneq Q^{-1}XQ^{-*}$, and partition $\widetilde{X}$ into blocks in the following two ways:
    \begin{equation}\label{eq:Xtilde_partition}
        \widetilde{X}=
        \begin{pmatrix}
            \widetilde{X}_{1,1} & \cdots & \widetilde{X}_{1,m} \\
            \vdots & \ddots & \vdots \\
            \widetilde{X}_{m,1} & \cdots & \widetilde{X}_{m,m}
        \end{pmatrix} \\
        =
        \begin{pmatrix}
            \widetilde{X}_{-,-} & \widetilde{X}_{-,0} & \widetilde{X}_{-,+} \\
            \widetilde{X}_{0,-} & \widetilde{X}_{0,0} & \widetilde{X}_{0,+} \\
            \widetilde{X}_{+,-} & \widetilde{X}_{+,0} & \widetilde{X}_{+,+}
        \end{pmatrix}.
    \end{equation}
    Furthermore, we denote the $(n_i,n_j)$-entry of $X_{i,j}$ by $\widetilde{x}_{i,j}$, which is the lower-right entry, and for $\ell=1,\ldots,m_0'$, let
    \begin{equation*}
        \widetilde{X}_{0,0}^{(\ell)}\coloneq
        \begin{pmatrix}
            \widetilde{x}_{\alpha+1,\alpha+1} & \cdots & \widetilde{x}_{\alpha+1,\beta} \\
            \vdots & \ddots & \vdots \\
            \widetilde{x}_{\beta,\alpha+1} & \cdots & \widetilde{x}_{\beta,\beta}
        \end{pmatrix}
        \in\CC^{m_0^{(\ell)}\times m_0^{(\ell)}},
    \end{equation*}
    where we set $\alpha=m_-+m_0^{(1)}+\cdots+m_0^{(\ell-1)},\ \beta=\alpha+m_0^{(\ell)}$ to simplify the indexing.

    Then, under the assumption that $W_\infty(x)=O$, it follows that $\widetilde{X}_{-,-}=O, \widetilde{X}_{+,+}=O$, and $\widetilde{X}_{0,0}^{(\ell)}=O \ (\ell=1,\ldots,m_0')$.
    Otherwise, the argument in \cref{subsec:scaled_gramian} would yield $W_\infty(x)\neq O$.
    Furthermore, we obtain
    \begin{gather*}
        \widetilde{X}_{-,-}=\sum_{i=1}^n x_i q_{-,i}q_{-,i}^*, \quad
        \widetilde{X}_{+,+}=\sum_{i=1}^n x_i q_{+,i}q_{+,i}^*, \\
        \widetilde{X}_{0,0}^{(\ell)}=\sum_{i=1}^n x_i q_{0,i}^{(\ell)}q_{0,i}^{(\ell)}{}^* \quad (\ell=1,\ldots,m_0').
    \end{gather*}
    By vectorizing the expression, we obtain $M_{\mathrm{VCS}}x=0$.
    Since $\rank M_{\mathrm{VCS}}=n$, it follows that $x=0$, which completes the proof.
    \qed
\end{proof}

\begin{proof}[\cref{prop:infinite_aecs_uniqueness}]
    Let $x=(x_1,\ldots,x_n)^\top\in\RR^n$ be any vector such that $W_-(x)=0$.
    Then, from \cref{prop:aecs_linearly_independent}, it suffices to prove $x=0$.
    Let $\widetilde{X}\coloneq Q^{-1}XQ^{-*}$ and partition $\widetilde{X}$ in the same way as \cref{eq:Xtilde_partition}.
    By an argument analogous to that of \cref{prop:infinite_vcs_uniqueness}, we obtain $\widetilde{X}_{-,-}=O$.
    Proceeding as in the proof of \cref{prop:infinite_vcs_uniqueness}, the assumption then implies that $x=0$, which completes the proof.
    \qed
\end{proof}

\section{Uniqueness of finite-horizon controllability scores}
\label{sec:finite_uniqueness}

In \cite{SatoTerasaki2024}, an example is presented in which the optimal solutions to \cref{prob:vcs,prob:aecs} are not unique for some $T>0$.
As noted in \cref{subsec:infinite_uniqueness}, the uniqueness of the optimal solutions is essential when interpreting controllability scores as network centrality.
In this section, we summarize existing results on the uniqueness of the controllability score on a finite time horizon.
We further present new results that encompass the existing ones and are applied to the Laplacian dynamics.

For Hurwitz matrices, it has been shown that the optimal solution is unique for all $T>0$.

\begin{proposition}[{\cite[Theorem 4]{SatoTerasaki2024}}]\label{prop:uniqueness_SatoTerasaki2024}
    Assume that $A$ is Hurwitz.
    Then, \cref{prob:vcs,prob:aecs} each admit a unique optimal solution for all $T>0$.
\end{proposition}

An important class of examples that does not satisfy this assumption is given by graph Laplacians.
Since a graph Laplacian has zero as an eigenvalue, it fails to satisfy the assumption.
For graph Laplacians of undirected connected graphs, the following result has been established.

\begin{proposition}[{\cite[Theorem 5]{SatoTerasaki2024}}]\label{prop:uniqueness_undirected}
    Let $L$ be the Laplacian matrix of an undirected connected graph, and let $A=-L$.
    Then, \cref{prob:vcs,prob:aecs} each admit a unique optimal solution for all $T>0$.
\end{proposition}

Without restricting the class of systems, the following result has been established for arbitrary matrices.

\begin{proposition}[{\cite[Theorem 1]{SatoKawamura2025}}]\label{prop:uniqueness_SatoKawamura2025}
    Let $A$ be an arbitrary matrix.
    Then, \cref{prob:vcs,prob:aecs} each admit a unique optimal solution for almost all $T>0$.
\end{proposition}

We present the following result, which encompasses all of the above results.

\begin{theorem}\label{prop:uniqueness_by_eigenvalue}
    Let $\lambda_i \ (i=1,\ldots,m)$ be the eigenvalues of $A$ as in \cref{subsec:settings},
    \begin{equation*}
        \Theta\coloneq
        \left\{
            \theta\in\RR\setminus\{0\}
        \,\middle|\,
            \text{$\lambda_i+\lambda_j=\theta\sqrt{-1}$ for some $i, j$}
        \right\},
    \end{equation*}
    and
    \begin{equation*}
        \Theta'\coloneq
        \left\{
            \left|\dfrac{2\ell\pi}{\theta}\right|\in\RR_{>0}
        \,\middle|\,
            \theta\in\Theta, \ell\in\ZZ\setminus\{0\}
        \right\}.
    \end{equation*}
    Then, \cref{prob:vcs,prob:aecs} each admit a unique optimal solution for $T\not\in\Theta'$.
\end{theorem}

To prove \cref{prop:uniqueness_by_eigenvalue}, we introduce some notation and prepare several lemmas.
For $Y=(y_1,\ldots,y_{\ell'}), \ y_i\in\RR^\ell$, we define
\begin{equation*}
    \vect(Y)\coloneq(y_1^\top,\ldots,y_{\ell'}^\top)^\top\in\RR^{\ell\ell'}.
\end{equation*}
Then, it is known that the following equality holds~\cite{Higham2008}:
\begin{equation}\label{eq:vec_identity}
    \vect(CYC'^\top)=(C'\otimes C)\vect(Y).
\end{equation}

Let $K\coloneq I_n\otimes A+A\otimes I_n\in\RR^{n^2\times n^2}$.
The following result is known for the spectrum of $K$ (see, e.g., \cite[Appendix B.14]{Higham2008}).
\begin{lemma}\label{prop:kronecker_spectrum}
    Every eigenvalue of $K$ is of the form $\lambda_i+\lambda_j$ for some $i,j$. 
\end{lemma}

We employ the following spectral mapping theorem.
\begin{lemma}[{\cite[Theorem 1.13]{Higham2008}}]\label{prop:spectral_mapping}
    Let $f\colon\CC\to\CC$ be a holomorphic function.
    Then, every eigenvalue of $f(K)$ is of the form $f(\mu)$, where $\mu$ is an eigenvalue of $K$.
\end{lemma}

Moreover, we employ the following lemma, which can be proved in the same way as in \cite[Theorems 1 and 3]{SatoTerasaki2024}.
\begin{lemma}\label{prop:linear_independence_lemma}
    Let $T>0$ be arbitrary.
    If $\widetilde{W}(x;T)=0$ for $x\in\RR^n$ implies $x=0$, then, \cref{prob:vcs,prob:aecs} each admit a unique optimal solution for that $T$.
\end{lemma}

\begin{proof}[\cref{prop:uniqueness_by_eigenvalue}]
    Let us consider a linear map
    \begin{equation}\label{eq:FT_definition}
        F_T\colon\RR^{n\times n}\to\RR^{n\times n}, \quad B\mapsto \widetilde{W}(B;T).
    \end{equation}

    Assume that $F_T(B)=0$ implies $B=0$.
    In particular, if $\widetilde{W}(x;T)=0$, then $x=0$, and by \cref{prop:linear_independence_lemma}, the optimal solution is unique for that $T$.
    We therefore investigate whether $F_T$ has zero as an eigenvalue.
    
    It follows from \cref{eq:vec_identity} that
    \begin{multline*}
        \vect\left(\widetilde{W}(B;T)\right)=\int_0^T \vect\left(\e^{At}B\e^{A^\top t}\right)\rd t \\
        =\left(\int_0^T \e^{At}\otimes \e^{At}\rd t\right)\vect(B)=\left(\int_0^T \e^{Kt}\rd t\right)\vect(B).
    \end{multline*}
    Thus, by vectorization, the map $F_T$ can be identified with the linear map
    \begin{equation}\label{eq:FT_identification}
        \RR^{n^2}\to\RR^{n^2}, \quad \vect(B)\mapsto\left(\int_0^T \e^{Kt}\rd t\right)\vect(B).
    \end{equation}

    Let
    \begin{equation*}
        f_T\colon \CC\to\CC, \quad z\mapsto\int_0^T \e^{zt}\rd t=
        \begin{cases*}
            \dfrac{\e^{zT}-1}{z} & if $z\neq 0$, \\
            T & if $z=0$,
        \end{cases*}
    \end{equation*}
    which is a holomorphic function.
    By identifying \cref{eq:FT_definition} and \cref{eq:FT_identification}, we have $F_T=f_T(K)$.
    Consequently, from \cref{prop:kronecker_spectrum,prop:spectral_mapping}, every eigenvalue of $F_T$ is of the form $f_T(\lambda_i+\lambda_j)$ for some $i, j$.
    
    Assume $f_T(\lambda_i+\lambda_j)=0$ for some $i, j$.
    Then, $\lambda_i+\lambda_j=\dfrac{2\ell\pi\sqrt{-1}}{T}$ for some $\ell\in\ZZ\setminus\{0\}$.
    In this case, $\lambda_i+\lambda_j$ is purely imaginary, and we can write $\lambda_i+\lambda_j=\theta\sqrt{-1}$ for some $\theta\neq0$.
    Consequently, $T=\dfrac{2\ell\pi}{\theta}=\left|\dfrac{2\ell\pi}{\theta}\right|$, and $T\in\Theta'$.

    Therefore, if $T\not\in\Theta'$, then $F_T$ does not have zero as an eigenvalue.
    It follows from \cref{prop:linear_independence_lemma} that \cref{prob:vcs,prob:aecs} each admit a unique optimal solution for this $T$.
    This completes the proof.
    \qed
\end{proof}

It follows that values of $T$ for which the optimal solution is not unique are confined to $\Theta'$.
Since $\Theta'$ is at most countable, the optimal solution is unique for almost all $T>0$.
In other words, \cref{prop:uniqueness_by_eigenvalue} can be understood as characterizing the exceptional values of $T>0$ for which the optimal solution fails to be unique in the result established by \cref{prop:uniqueness_SatoKawamura2025}.

Moreover, although \cref{prop:uniqueness_undirected} is stated for Laplacians of undirected connected graphs, the same uniqueness property is shown to hold for directed graphs as well.

\begin{corollary}\label{prop:uniqueness_directed}
    Let $L$ be the Laplacian matrix of a directed graph, and let $A=-L$.
    Then, \cref{prob:vcs,prob:aecs} each admit a unique optimal solution for all $T>0$.
\end{corollary}

\begin{proof}
    It follows from \cref{prop:laplacian_eigenvalue} that $\Theta=\emptyset$.
    Therefore, $\Theta'=\emptyset$ also holds, and this completes the proof.
    \qed
\end{proof}

\section*{Acknowledgment}
This work was supported by JST PRESTO, Japan, Grant Number JPMJPR25K4.

\bibliographystyle{IEEEtran}
\bibliography{reference.bib}

@Misc{amsmath,
  author =	 {{American Mathematical Society}},
  title =	 {User's Guide for the \texttt{amsmath} Package
                  (Version 2.0)},
  url =		 {ftp://ftp.ams.org/pub/tex/doc/amsmath/amsldoc.pdf},
  urldate =	 {2015-07-30},
  year =	 2002}

@Misc{pgfplots,
  author =	 {Christian Feuers\"anger},
  title =	 {Manual for Package \texttt{PGFPLOTS}},
  month =	 may,
  year =	 2015,
  url =		 {http://sourceforge.net/projects/pgfplots}
}

@article {SatoTerasaki2024,
    AUTHOR = {Sato, Kazuhiro and Terasaki, Shun},
     TITLE = {Controllability scores for selecting control nodes of
              large-scale network systems},
   JOURNAL = {IEEE Trans. Automat. Control},
  FJOURNAL = {Institute of Electrical and Electronics Engineers.
              Transactions on Automatic Control},
    VOLUME = {69},
      YEAR = {2024},
    NUMBER = {7},
     PAGES = {4673--4680},
      ISSN = {0018-9286,1558-2523},
   MRCLASS = {93B51 (90C25 93A15 93B05)},
  MRNUMBER = {4765743},
       DOI = {10.1109/tac.2024.3355806},
       URL = {https://doi.org/10.1109/tac.2024.3355806},
}

@article {SatoKawamura2025,
    AUTHOR = {Sato, Kazuhiro and Kawamura, Ryohei},
     TITLE = {Uniqueness {A}nalysis of {C}ontrollability {S}cores and
              {T}heir {A}pplication to {B}rain {N}etworks},
   JOURNAL = {IEEE Trans. Control Netw. Syst.},
  FJOURNAL = {IEEE Transactions on Control of Network Systems},
    VOLUME = {12},
      YEAR = {2025},
    NUMBER = {4},
     PAGES = {2568--2580},
      ISSN = {2325-5870},
   MRCLASS = {93B05 (05 92 93C05)},
  MRNUMBER = {5007236},
       DOI = {10.1109/tcns.2025.3583613},
       URL = {https://doi.org/10.1109/tcns.2025.3583613},
}

@article {Liu2011,
    AUTHOR = {Liu, Yang-Yu and Slotine, Jean-Jacques and Barab\'{a}si, Albert-L\'{a}szl\'{o}},
     TITLE = {Controllability of complex networks},
   JOURNAL = {Nature},
  FJOURNAL = {Nature},
    VOLUME = {473},
      YEAR = {2011},
     PAGES = {167--173},
      ISSN = {0028-0836,1476-4687},
       DOI = {doi.org/10.1038/nature10011},
       URL = {https://doi.org/10.1038/nature10011},
}

@article {Lin1974,
    AUTHOR = {Lin, Ching Tai},
     TITLE = {Structural controllability},
   JOURNAL = {IEEE Trans. Automatic Control},
  FJOURNAL = {Institute of Electrical and Electronics Engineers.
              Transactions on Automatic Control},
    VOLUME = {AC-19},
      YEAR = {1974},
     PAGES = {201--208},
   MRCLASS = {93B05},
  MRNUMBER = {452870},
       DOI = {10.1109/tac.1974.1100557},
       URL = {https://doi.org/10.1109/tac.1974.1100557},
}

@article {Summers2016,
    AUTHOR = {Summers, Tyler H. and Cortesi, Fabrizio L. and Lygeros, John},
     TITLE = {On submodularity and controllability in complex dynamical
              networks},
   JOURNAL = {IEEE Trans. Control Netw. Syst.},
  FJOURNAL = {IEEE Transactions on Control of Network Systems},
    VOLUME = {3},
      YEAR = {2016},
    NUMBER = {1},
     PAGES = {91--101},
      ISSN = {2325-5870},
   MRCLASS = {93B05 (90C27 93B07 93B51 93C05)},
  MRNUMBER = {3477965},
MRREVIEWER = {M.\ Amin\ Rahimian},
       DOI = {10.1109/TCNS.2015.2453711},
       URL = {https://doi.org/10.1109/TCNS.2015.2453711},
}

@article {Gu2015,
    AUTHOR = {Gu, Shi and Pasqualetti, Fabio and Cieslak, Matthew and Telesford, Qawi K. and Yu, Alfred B. and Kahn, Ari E. and Medaglia, John D. and Vettel, Jean M. and Miller, Michael B. and Grafton, Scott T. and Bassett, Danielle S.},
     TITLE = {Controllability of structural brain networks},
   JOURNAL = {Nat. Commun.},
  FJOURNAL = {Nature Communications},
    VOLUME = {6},
      YEAR = {2015},
    NUMBER = {1},
       EID = {8414},
      ISSN = {2041-1723},
       DOI = {10.1038/ncomms9414},
       URL = {https://doi.org/10.1038/ncomms9414},
 PAGETOTAL = {10},
}

@article{Liu2012,
    AUTHOR = {Liu, Yang-Yu and Slotine, Jean-Jacques and Barab\'{a}si, Albert-L\'{a}szl\'{o}},
     TITLE = {Control centrality and hierarchical structure in complex networks},
   JOURNAL = {PLoS One},
  FJOURNAL = {PLoS One},
    VOLUME = {7},
      YEAR = {2012},
    NUMBER = {9},
       EID = {e44459},
      ISSN = {1932-6203},
       DOI = {10.1371/journal.pone.0044459},
       URL = {https://doi.org/10.1371/journal.pone.0044459},
 PAGETOTAL = {7},
}

@book {HornJohnson2012,
    AUTHOR = {Horn, Roger A. and Johnson, Charles R.},
     TITLE = {Matrix Analysis},
   EDITION = {2},
 PUBLISHER = {Cambridge University Press},
     PLACE = {Cambridge},
      YEAR = {2012},
     PAGES = {xviii+662},
      ISBN = {978-0-521-83940-2},
}

@book {RockafellarWets1998,
    AUTHOR = {Rockafellar, R. Tyrrell and Wets, Roger J.-B.},
     TITLE = {Variational {A}nalysis},
    SERIES = {Grundlehren der mathematischen Wissenschaften [Fundamental
              Principles of Mathematical Sciences]},
    VOLUME = {317},
 PUBLISHER = {Springer-Verlag, Berlin},
      YEAR = {1998},
     PAGES = {xiv+733},
      ISBN = {3-540-62772-3},
   MRCLASS = {49-02 (46N10 47N10 49J52 49K40 90C30)},
  MRNUMBER = {1491362},
MRREVIEWER = {Francis\ H.\ Clarke},
       DOI = {10.1007/978-3-642-02431-3},
       URL = {https://doi.org/10.1007/978-3-642-02431-3},
}

@article {MullerWeber1972,
    AUTHOR = {M\"uller, P. C. and Weber, H. I.},
     TITLE = {Analysis and optimization of certain qualities of
              controllability and observability for linear dynamical
              systems},
   JOURNAL = {Automatica},
  FJOURNAL = {Automatica},
    VOLUME = {8},
      YEAR = {1972},
    NUMBER = {3},
     PAGES = {237--246},
      ISSN = {0005-1098,1873-2836},
   MRCLASS = {93B05},
  MRNUMBER = {335016},
MRREVIEWER = {T.\ Kaczorek},
       DOI = {10.1016/0005-1098(72)90044-1},
       URL = {https://doi.org/10.1016/0005-1098(72)90044-1},
}

@article{SatoOpenDomain2026,
  author  = {Kazuhiro Sato},
  title   = {Projected-Gradient Analysis for Open-Domain Convex
             Optimization under Boundary Blow-Up: Application to
             Controllability Scoring},
  journal = {arXiv preprint arXiv:2607.23674},
  year    = {2026}
}

@article{umezu2026controllability,
  title={Controllability scores of linear time-varying network systems},
  author={Umezu, Kota and Sato, Kazuhiro},
  journal={IEEE Transactions on Automatic Control},
   volume  = {71},
  number  = {8},
  pages   = {5144--5159},
  year={2026},
  publisher={IEEE},
  doi = {10.1109/TAC.2026.3669811}
}

@misc{ControllabilityScoreCode,
  author       = {Sato, Kazuhiro and Umezu, Kota and Yamada, Ritsuki},
  title        = {{Controllability Score Computation Program}},
  howpublished = {GitHub repository},
  year         = {2026},
  note         = {[Online]. Available:
                  \url{https://github.com/Sato-Sakigake-Controllability-Score/controllability-score}.
                  Accessed: Aug. 6, 2026}
}

@article{Simoncini2007,
  author  = {Valeria Simoncini},
  title   = {A New Iterative Method for Solving Large-Scale {Lyapunov} Matrix Equations},
  journal = {SIAM Journal on Scientific Computing},
  volume  = {29},
  number  = {3},
  pages   = {1268--1288},
  year    = {2007},
  doi     = {10.1137/06066120X}
}

@article{BennerKurschnerSaak2014,
  author  = {Peter Benner and Patrick K{\"u}rschner and Jens Saak},
  title   = {Self-Generating and Efficient Shift Parameters in {ADI} Methods for Large {Lyapunov} and {Sylvester} Equations},
  journal = {Electronic Transactions on Numerical Analysis},
  volume  = {43},
  pages   = {142--162},
  year    = {2014}
}

@inproceedings{SatoOED2026,
  author    = {Kazuhiro Sato},
  title     = {Relationship Between Controllability Scoring and Optimal Experimental Design},
  booktitle = {Proceedings of the 15th Asian Control Conference (ASCC)},
  address   = {Bali, Indonesia},
  year      = {2026},
  note      = {Available at arXiv:2602.11921},
  url       = {https://arxiv.org/abs/2602.11921}
}

@article {ChebotarevAgaev2002,
    AUTHOR = {Chebotarev, Pavel and Agaev, Rafig},
     TITLE = {Forest matrices around the {L}aplacian matrix},
   JOURNAL = {Linear Algebra Appl.},
  FJOURNAL = {Linear Algebra and its Applications},
    VOLUME = {356},
      YEAR = {2002},
     PAGES = {253--274},
      ISSN = {0024-3795,1873-1856},
   MRCLASS = {05C05 (05C20)},
  MRNUMBER = {1945250},
MRREVIEWER = {W.-K.\ Chen},
       URL ={https://www.sciencedirect.com/science/article/pii/S0024379502003889},
}

@book {Varga2004,
    AUTHOR = {Varga, Richard S.},
     TITLE = {Ger\v sgorin and his circles},
    SERIES = {Springer Series in Computational Mathematics},
    VOLUME = {36},
 PUBLISHER = {Springer-Verlag, Berlin},
      YEAR = {2004},
     PAGES = {x+226},
      ISBN = {3-540-21100-4},
   MRCLASS = {15-02 (15-03 15A18 15A42 65-02 65-03 65F15)},
  MRNUMBER = {2093409},
MRREVIEWER = {David\ Scott\ Watkins},
       DOI = {10.1007/978-3-642-17798-9},
       URL = {https://doi.org/10.1007/978-3-642-17798-9},
}

@article {Lancaster1970,
    AUTHOR = {Lancaster, Peter},
     TITLE = {Explicit solutions of linear matrix equations},
   JOURNAL = {SIAM Rev.},
  FJOURNAL = {SIAM Review. A Publication of the Society for Industrial and
              Applied Mathematics},
    VOLUME = {12},
      YEAR = {1970},
     PAGES = {544--566},
      ISSN = {1095-7200},
   MRCLASS = {15.35},
  MRNUMBER = {279115},
MRREVIEWER = {A.\ S.\ Householder},
       DOI = {10.1137/1012104},
       URL = {https://doi.org/10.1137/1012104},
}

@article{belabbas2025infinite,
  title={On Infinite-horizon Minimum Energy Control},
  author={Belabbas, Mohamed-Ali and Chen, Xudong},
  journal={arXiv preprint arXiv:2502.17795},
  year={2025}
}

@article{Sato2025,
  title={Target Controllability Scores for Actuation-Constrained Network Intervention},
  author={Sato, Kazuhiro},
  journal={arXiv preprint arXiv:2510.13354},
  year={2026}
}

@book {Higham2008,
    AUTHOR = {Higham, Nicholas J.},
     TITLE = {Functions of matrices},
 PUBLISHER = {Society for Industrial and Applied Mathematics (SIAM),
              Philadelphia, PA},
      YEAR = {2008},
     PAGES = {xx+425},
      ISBN = {978-0-89871-646-7},
   MRCLASS = {15-02 (93B40 93C05)},
  MRNUMBER = {2396439},
MRREVIEWER = {Daniel\ Kressner},
       DOI = {10.1137/1.9780898717778},
       URL = {https://doi.org/10.1137/1.9780898717778},
}

@inproceedings {Kalman1960,
    AUTHOR = {Kalman, R. E.},
 BOOKTITLE = {Proceedings First International Conference on Automatic Control},
     TITLE = {On the general theory of control systems},
      YEAR = {1960},
     PAGES = {481--492},
}

@article {Newman2003,
    AUTHOR = {Newman, M. E. J.},
     TITLE = {The structure and function of complex networks},
   JOURNAL = {SIAM Rev.},
  FJOURNAL = {SIAM Review},
    VOLUME = {45},
      YEAR = {2003},
    NUMBER = {2},
     PAGES = {167--256},
      ISSN = {1095-7200,0036-1445},
   MRCLASS = {05C90 (68M10 90B10 90B18 91D30 92B20 94C15)},
  MRNUMBER = {2010377},
MRREVIEWER = {Jerrold\ W.\ Grossman},
       DOI = {10.1137/S003614450342480},
       URL = {https://doi-org/10.1137/S003614450342480},
}

@article {Yan2015,
    AUTHOR = {Yan, Gang and Tsekenis, Georgios and Barzel, Baruch and Slotine, Jean-Jacques and Liu, Yang-Yu and Barab{\'a}si, Albert-L{\'a}szl{\'o}},
     TITLE = {Spectrum of controlling and observing complex networks},
   JOURNAL = {Nat. Phys.},
  FJOURNAL = {Nature Physics},
    VOLUME = {11},
      YEAR = {2015},
    NUMBER = {9},
     PAGES = {779--786},
      ISSN = {1745-2481,1745-2473},
       DOI = {10.1038/nphys3422},
       URL = {https://doi.org/10.1038/nphys3422},
}

@article {BresslerMenon2010,
    AUTHOR = {Bressler, Steven L. and Menon, Vinod},
     TITLE = {Large-scale brain networks in cognition: emerging methods and principles},
   JOURNAL = {Trends Cogn. Sci.},
  FJOURNAL = {Trends in Cognitive Sciences},
    VOLUME = {14},
      YEAR = {2010},
    NUMBER = {6},
     PAGES = {277--290},
      ISSN = {1364-6613},
       DOI = {10.1016/j.tics.2010.04.004},
       URL = {https://doi.org/10.1016/j.tics.2010.04.004},
}

@article {WigSchlaggarPetersen2011,
    AUTHOR = {Wig, Gagan S. and Schlaggar, Bradley L. and Petersen, Steven E.},
     TITLE = {Concepts and principles in the analysis of brain networks},
   JOURNAL = {Ann. N. Y. Acad. Sci.},
  FJOURNAL = {Annals of the New York Academy of Sciences},
    VOLUME = {1224},
      YEAR = {2011},
    NUMBER = {1},
     PAGES = {126--146},
      ISSN = {1749-6632,0077-8923},
       DOI = {10.1111/j.1749-6632.2010.05947.x},
       URL = {https://doi.org/10.1111/j.1749-6632.2010.05947.x},
}

@article {BarabasiGulbahceLoscalzo2011,
    AUTHOR = {Barab{\'a}si, Albert-L{\'a}szl{\'o} and Gulbahce, Natali and Loscalzo, Joseph},
     TITLE = {Network medicine: a network-based approach to human disease},
   JOURNAL = {Nat. Rev. Genet.},
  FJOURNAL = {Nature Reviews Genetics},
    VOLUME = {12},
      YEAR = {2011},
    NUMBER = {1},
     PAGES = {56--68},
      ISSN = {1471-0064},
       DOI = {10.1038/nrg2918},
       URL = {https://doi.org/10.1038/nrg2918},
}

@inproceedings{Backstrom2012,
    AUTHOR = {Backstrom, Lars and Boldi, Paolo and Rosa, Marco and Ugander, Johan and Vigna, Sebastiano},
     TITLE = {Four degrees of separation},
 BOOKTITLE = {Proceedings of the 4th Annual ACM Web Science Conference},
      YEAR = {2012},
     PAGES = {33--42},
       DOI = {10.1145/2380718.2380723},
       URL = {https://doi.org/10.1145/2380718.2380723},
}

@article {Buldyrev2010,
    AUTHOR = {Buldyrev, Sergey V. and Parshani, Roni and Paul, Gerald and Stanley, H. Eugene and Havlin, Shlomo},
     TITLE = {Catastrophic cascade of failures in interdependent networks},
   JOURNAL = {Nature},
  FJOURNAL = {Nature},
    VOLUME = {464},
      YEAR = {2010},
    NUMBER = {7291},
     PAGES = {1025--1028},
      ISSN = {1476-4687,0028-0836},
       DOI = {10.1038/nature08932},
       URL = {https://doi.org/10.1038/nature08932},
}

@article {PaganiAiello2013,
    AUTHOR = {Pagani, Giuliano Andrea and Aiello, Marco},
     TITLE = {The Power Grid as a complex network: A survey},
   JOURNAL = {Phys. A: Stat. Mech. Appl.},
  FJOURNAL = {Physica A: Statistical Mechanics and its Applications},
    VOLUME = {392},
      YEAR = {2013},
    NUMBER = {11},
     PAGES = {2688--2700},
      ISSN = {0378-4371},
       DOI = {10.1016/j.physa.2013.01.023},
       URL = {https://doi.org/10.1016/j.physa.2013.01.023},
}

@article {ElterenQuaxSloot2022,
    AUTHOR = {van Elteren, Casper and Quax, Rick and Sloot, Peter},
     TITLE = {Dynamic importance of network nodes is poorly predicted by static structural features},
   JOURNAL = {Phys. A: Stat. Mech. Appl.},
  FJOURNAL = {Physica A: Statistical Mechanics and its Applications},
    VOLUME = {593},
      YEAR = {2022},
    eprint = {126889},
      ISSN = {0378-4371},
       DOI = {10.1016/j.physa.2022.126889},
       URL = {https://doi.org/10.1016/j.physa.2022.126889},
}

@article {Vinayagam2016,
    AUTHOR = {Vinayagam, Arunachalam and Gibson, Travis E. and Lee, Ho-Joon and Yilmazel, Bahar and Roesel, Charles and Hu, Yanhui and Kwon, Young and Sharma, Amitabh and Liu, Yang-Yu and Perrimon, Norbert and Barab{\'a}si, Albert-L{\'a}szl{\'o}},
     TITLE = {Controllability analysis of the directed human protein interaction network identifies disease genes and drug targets},
   JOURNAL = {Proc. Natl. Acad. Sci. U.S.A.},
  FJOURNAL = {Proceedings of the National Academy of Sciences of the United States of America},
    VOLUME = {113},
      YEAR = {2016},
    NUMBER = {18},
     PAGES = {4976--4981},
      ISSN = {0378-4371},
       DOI = {10.1073/pnas.1603992113},
       URL = {https://doi.org/10.1073/pnas.1603992113},
}

@article {GolubWilkinson1976,
    AUTHOR = {Golub, G. H. and Wilkinson, J. H.},
     TITLE = {Ill-conditioned eigensystems and the computation of the
              {J}ordan canonical form},
   JOURNAL = {SIAM Rev.},
  FJOURNAL = {SIAM Review. A Publication of the Society for Industrial and
              Applied Mathematics},
    VOLUME = {18},
      YEAR = {1976},
    NUMBER = {4},
     PAGES = {578--619},
      ISSN = {1095-7200},
   MRCLASS = {65F15 (15A21)},
  MRNUMBER = {413456},
MRREVIEWER = {Robert\ Todd\ Gregory},
       DOI = {10.1137/1018113},
       URL = {https://doi-org/10.1137/1018113},
}

\end{document}